\documentclass[letter,onetabnum]{siamart220329}

\usepackage{booktabs}
\usepackage[disable,obeyFinal,textsize=footnotesize]{todonotes}

\newcommand{\ekin}[2][]{\todo[color=green, #1]{\textbf{Ekin}: #2}}

\usepackage{textcomp}

\usepackage{amsmath}
\newtheorem{defn}{Definition}[section]

\usepackage{array}
\usepackage{fancyvrb}
\usepackage{algpseudocode}
\usepackage{algorithm}
\usepackage{listings}
\usepackage{mdframed}
\usepackage{graphicx,grffile}
\usepackage{amsmath,amssymb}
\usepackage{resizegather}
\usepackage{tikz}
\usepackage{subcaption}

\definecolor{ekinblue}{RGB}{75,172,198}
\definecolor{ekinorange}{RGB}{247,150,70}
\definecolor{ekinpurple}{RGB}{112,48,160}
\usepackage{xcolor}

 \makeatletter
\def\l@subsection{\@tocline{2}{0pt}{2.5pc}{5pc}{}}
\def\l@subsubsection{\@tocline{2}{0pt}{3.5pc}{5pc}{}}
 \makeatletter

\makeatletter
\newcommand\notsotiny{\@setfontsize\notsotiny\@vipt\@viipt}
\makeatother

\begin{document}

\title{Backpropagation through Back Substitution with a Backslash}

\author{Alan Edelman\thanks{Department of Mathematics and CSAIL, MIT, Cambridge, MA
    (\email{edelman@mit.edu}).}
  \and Ekin Aky\"urek\thanks{Department of EECS and CSAIL, MIT, Cambridge, MA  \email{(akyurek@mit.edu)}}
  \and Yuyang Wang\thanks{AWS AI Labs, Santa Clara, CA
    (\email{yuyawang@amazon.com}). Work done prior to joining Amazon.}}

\maketitle

\begin{abstract} We present a linear algebra formulation of backpropagation which allows the calculation of gradients by using a generically written ``backslash'' or Gaussian elimination on triangular systems of equations.  Generally, the matrix elements are  operators. This paper has three contributions: (i) it is of intellectual value to replace traditional treatments of automatic differentiation with a (left acting) operator theoretic, graph-based approach; (ii) operators can be readily placed in matrices in software in programming languages such as Julia as an implementation option; (iii) we introduce a novel notation, ``transpose dot'' operator ``$\{\}^{T_\bullet}$'' that allows for the reversal of operators.

We further demonstrate  the elegance of the operators approach in a suitable programming language consisting of generic linear algebra operators such as Julia \cite{bezanson2017julia}, and that it is possible to realize this abstraction in code. Our implementation shows how generic linear algebra can allow operators as elements of matrices. In contrast to ``operator overloading,'' where backslash  would normally have to be rewritten to take advantage of operators, with ``generic programming'' there is no such need.

\end{abstract}

\section{Preface: Summary and the Challenge}
\label{preface}

This paper provides the mathematics to show how an operator theoretic, graph-based approach can realize backpropagation by applying back substitution to a matrix whose elements are operators.

As a showcase result, one can back-propagate to compute the gradient on feed-forward neural networks (or Multi-layer Perceptron (MLP))~\cite{Goodfellow-et-al-2016} with 
\begin{equation}
  \label{eqn:bbb}
  \nabla J = M^T((I-\tilde{L})^T \backslash g),
\end{equation}
where $M$ (source to non-source nodes) and $\tilde{L}$ (within non-source nodes) are blocks of the adjacency matrix of the computational graph (see Section~\ref{sec:mlp} for precise definitions), $g$ is the vector of gradients of the loss function, and $I$ is the identity matrix. For readers unfamiliar with the backslash notation, an equivalent expression of \eqref{eqn:bbb} is $\nabla J = M^T(I-\tilde{L})^{-T} g$.

We then set up a challenge to ourselves. Could we correctly implement \eqref{eqn:bbb} by simply typing the command (after basic setup but \emph{without overloading of backslash})

\vspace{.1in}

\begin{center}
  \includegraphics[width=1.6in]{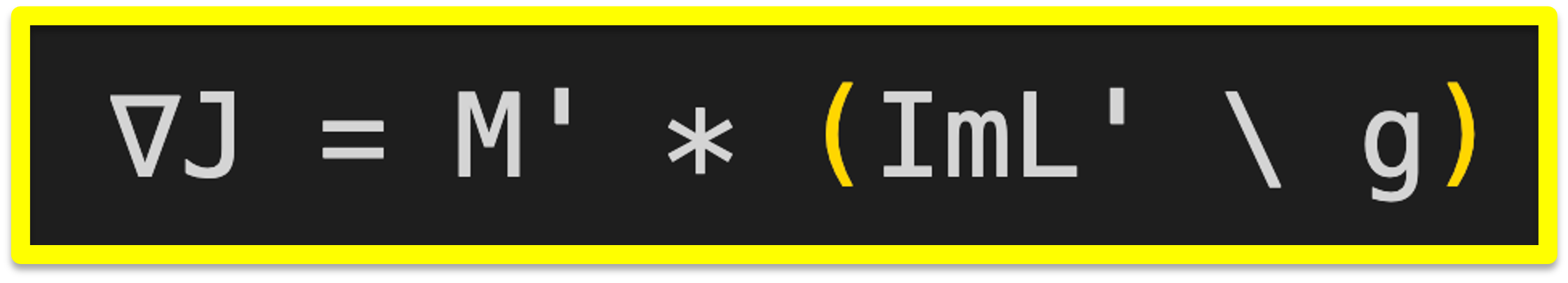} \raisebox{.1in}?
\end{center}

We demonstrate that indeed the backpropagation can be achieved, almost by magic, in a programming language armed with generic programming allowing for operators as elements of matrices. The software in
Section \ref{challenge} is by itself interesting not for the usual reasons of what it does, but in this case how it does it:
how a powerful language with generic programming and multiple dispatch can allow this abstract mathematical formulation to be realized.

\section{Introduction: Linear Algebra, Graphs, Automatic Differentiation (AD), Operators, and Julia}

Automatic differentiation (AD) is fundamental to gradient-based
optimization of neural networks and is used throughout scientific computing.
There are two popular approaches to AD: namely, forward and backward
(reverse) modes~\cite{griewank2003mathematical,griewank2008evaluating,JMLR:v18:17-468, revels2016forward}, the latter of which is also known as \emph{backpropagation} in the Machine Learning (ML) literature.\footnote{
  Despite similar terminology, the term ``forward propagation'' (or forward pass) in machine learning (ML) has no connection to forward mode automatic differentiation (AD). Instead, it refers to the process where a neural network calculates its output by sequentially passing input data through each layer, applying weighted sums and activation functions, until it reaches the output layer. In later sections, Algorithm~\ref{nnalg} and~\ref{matrixneuralnet} illustrate such a procedure. Whereas, ``backpropagation'' (backward pass, or reverse mode AD) is so named because information flows backwards through the network during this process.
}
A common high-level description of AD is that
it is really ``only''  the chain-rule. The centuries old technology of taking derivatives
is taking on a modern twist in the form of Differentiable Programming~\cite{Wik:DifProg,li2018differentiable}.
Who would have thought that one of the most routine college course subjects would now be the
subject of  much renewed interest both in applied mathematics and computer science?

This paper introduces the notion that AD is best understood with a matrix-based approach. The chain-rule explanation, in retrospect, feels to us as a distraction or at least extra baggage. We suspect that while the chain-rule is well known, it is understood mechanically rather than deeply by most of us. We argue that a linear algebra based framework for AD, while mathematically equivalent to other approaches, provides a simplicity of understanding, and equally importantly a viable approach worthy of further study.

Regarding software, while most high-level languages allow for matrices whose elements are scalars, the ability to work with matrices whose elements might be operators without major changes to the elementwise software is an intriguing abstraction. We discuss a Julia implementation that makes this step particularly mathematically natural.

It is our view that a linear algebraic approach sheds light on how backpropagation works in its essence.
We theoretically connect backpropagation to the back substitution method for triangular systems of equations. Similarly, forward substitution corresponds to the forward mode calculation of automatic differentiation.
As is well documented in the preface to the book Graph Algorithms in the Language of Linear Algebra \cite{kepner2011graph}, there have been many known benefits to formulate mathematically a graph algorithm in linear algebraic terms.

The ability to implement these abstractions while retaining performance is demonstrated using Julia, a language that facilitates abstractions, multiple dispatch, the type system, and which offers generic operators.

\section{A Matrix Method for Weighted Paths}

\subsection{``Forward and Back'' through Graphs and Linear Algebra}
\label{frg}

In the spirit of reaching the  mathematical core,
let us strip away the derivatives, gradients, Jacobians,  the computational graphs, and the ``chain rule''  that clutter the story of how is it possible to compute the same thing forwards and backwards.
We set ourselves the goal of explaining the essence of forward mode vs backward mode in AD with a single figure. Figure \ref{hw} is the result. Note that ``forward mode'' differentiation is not to be confused with the forward computation of the desired quantity.

\subsubsection{Path Weights on Directed Graphs} 
\begin{figure}[H]
  \begin{center}
    \includegraphics[width=.9\textwidth]{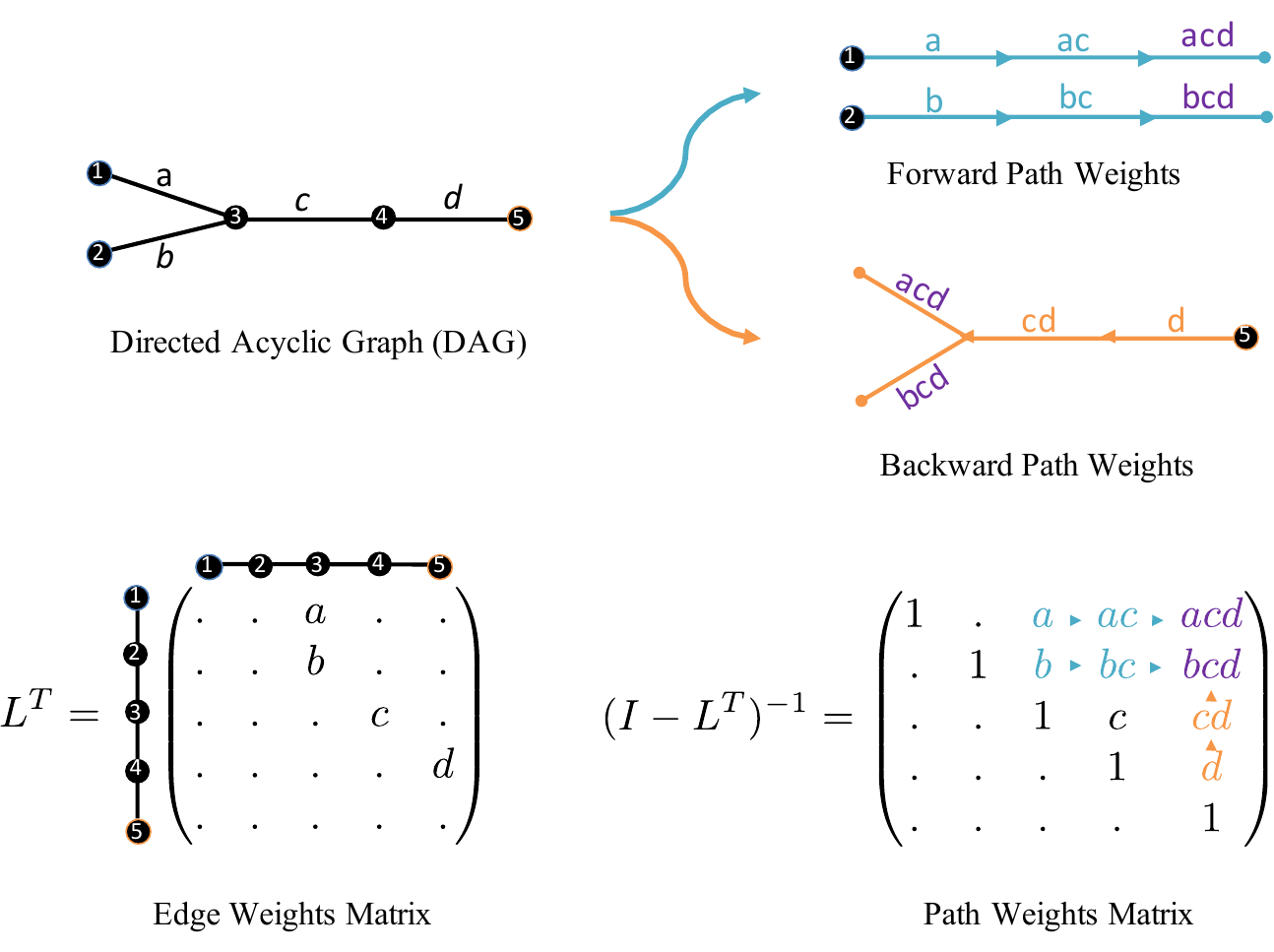}    
    \caption{
    {{\bf Legend} : {\color{ekinpurple} Purple: target weights}, {\color{ekinblue} Blue: forward computation}, {\color{ekinorange} Orange: backward computation}. The dots in matrices denote zeros.} \\
    (Upper Left:) Multiply the weights along the paths from source node 1  to sink node 5 and also source node 2 to sink node 5 to obtain
      {\color{ekinpurple}{acd}} and {\color{ekinpurple}{bcd}}.
    {\color{ekinblue} (Right Blue:) The obvious forward method.}
    {\color{ekinorange} (Right Orange:) A backward method that requires one fewer multiplication. }\\
    (Below:) A matrix method: if
    $L^T_{ij}=$
    the weight on  edge $ij$, then $(I-L^T)^{-1}$ simultaneously exhibits the {\color{ekinblue} forward} (i.e. $a\rightarrow ac \rightarrow acd$ and $b\rightarrow bc \rightarrow bcd$) and {\color{ekinorange} backward} methods (i.e. $d\rightarrow cd \rightarrow bcd \rightarrow acd / bcd$). 
    }
    \label{hw}
  \end{center}
\end{figure}

Consider a directed acyclic graph (DAG) with edge weights as
in Figure \ref{hw} where nodes 1 and 2 are sources (starting nodes), and node 5 is a sink (end node). The problem is to compute the {\bf path weights}, which we define as the products of the weights from every start node to every sink node.

Evidently, the path weights that we seek in Figure \ref{hw} may be obtained by calculating
\begin{equation}
  \label{assoc}
  \mbox{path weights} =
  \underbrace{\begin{pmatrix} 1 & 0  \\  0 & 1 \\ 0 & 0  \\  0 & 0 \\ 0 & 0 \end{pmatrix}^{\! T}}_{\mbox{sources}}
  (I-L^T)^{-1}
  \underbrace{\begin{pmatrix} 0 \\  0 \\ 0 \\  0 \\ 1 \end{pmatrix}}_{\mbox{sink}} ,
\end{equation}
where $L^T$, the adjacency matrix (edge weight matrix), is displayed in the lower left of Figure \ref{hw}.
One explanation of why
\eqref{assoc}
works for calculating the path weights is that $(L^T)^k_{ij}$ sums the path weights of length $k$ from node $i$ to node $j$
and  $(I-L^T)^{-1} = I + L^T + \ldots + (L^T)^{n-1}$\ekin{I guess this is true for this specific matrix, should we mention that?}  then counts  path weights of all lengths from $i$ to $j$.
\definecolor{dg}{rgb}{0,0.5,0}

If one follows step by step the linear algebra methods of forward substitution for lower triangular matrices or back substitution for upper triangular matrices, one obtains path weights algorithms as summarized in Figure \ref{forbac}. We remind the reader that forward and back substitution are the standard methods to solve lower and upper triangular systems respectively.

\begin{figure}[h]
  $$\begin{array}{|lc|}
      \hline
      \multicolumn{2}{|l|}{\mbox{Two Equivalent Ways to Compute the Path Weights in Figure \ref{hw}:}}                   \\
      \hline                        &                                                                                    \\[.03in]
      \mbox{Forward Substitution: } & \begin{pmatrix} 0 \\  0 \\ 0 \\  0 \\ 1 \end{pmatrix}^{\! \! \! \!  T} \!
      \ \left[ (I-L)^{-1} \begin{pmatrix} 1 & 0  \\  0 & 1 \\ 0 & 0  \\  0 & 0 \\ 0 & 0 \end{pmatrix}  \right]  \! \! \! \\[.5in]
      \mbox{Back  Substitution: }   &
      \begin{pmatrix} 1 & 0  \\  0 & 1 \\ 0 & 0  \\  0 & 0 \\ 0 & 0 \end{pmatrix}^{\! \! \! \!  T}     \! \!
      \left[ (I-L^T)^{-1}    \begin{pmatrix} 0 \\  0 \\ 0 \\  0 \\ 1 \end{pmatrix} \ \   \right]                         \\[.5in]   \hline
    \end{array}$$
  \caption{\label{forbac}The forward and backward methods
    compared: Both are seen equivalently as a choice of parenthesizing \eqref{assoc}
    or as forward substitution vs.\ back substitution.  Generally speaking, when the number of sources is larger than the number of sinks,
    one might expect the backward  method to have less complexity.}
\end{figure}

\subsubsection{Generalizing ``Forward and Back'' to a Catalan number of possibilities}

Continuing with the same $L$ matrix from Section \ref{frg},  we can
begin to understand all of the possibilities including the forward method,
the backward method, the mixed-modes methods, and even more possibilities:
\[
  \! \!  \! \! \!
  \begin{array}{c}
    \begin{pmatrix} 1 & 0  \\  0 & 1 \\ 0 & 0  \\  0 & 0 \\ 0 & 0 \end{pmatrix}^{\! T}  (I-L^T)^{-1} \begin{pmatrix} 0 \\  0 \\ 0 \\  0 \\ 1 \end{pmatrix} = \\ \\
    \begin{pmatrix} 1 & 0  \\  0 & 1 \\ 0 & 0  \\  0 & 0 \\ 0 & 0 \end{pmatrix}^{\! \!  \!  \!T} \! \!  \! \!
    \begin{pmatrix}
      1 & . & a & . & . \\
      . & 1 & . & . & . \\
      . & . & 1 & . & . \\
      . & . & . & 1 & . \\
      . & . & . & . & 1 \\
    \end{pmatrix}\! \!  \! \!
    \begin{pmatrix}
      1 & . & . & . & . \\
      . & 1 & b & . & . \\
      . & . & 1 & . & . \\
      . & . & . & 1 & . \\
      . & . & . & . & 1 \\
    \end{pmatrix} \! \!  \! \!
    \begin{pmatrix}
      1 & . & . & . & . \\
      . & 1 & . & . & . \\
      . & . & 1 & c & . \\
      . & . & . & 1 & . \\
      . & . & . & . & 1 \\
    \end{pmatrix}\! \!  \! \!
    \begin{pmatrix}
      1 & . & . & . & . \\
      . & 1 & . & . & . \\
      . & . & 1 & . & . \\
      . & . & . & 1 & d \\
      . & . & . & . & 1 \\
    \end{pmatrix}\! \!  \! \!
    \begin{pmatrix} 0 \\  0 \\ 0 \\  0 \\ 1 \end{pmatrix}
    \!  \!   .
  \end{array}
\]

It is well known \cite{stanley2015catalan}, that there are a Catalan number, $C_5=42$, ways to parenthesize the above expression.
One of the 42 choices evaluates left to right; this is forward substitution  which computes the graph weights forward.
Another evaluating from right to left is backward substitution.   There are three other ``mixed-modes''  \cite{1810.08297}
which combine forward and backward methods.  The remaining 37 methods require matrix-matrix multiplication as a first step.
We encourage the reader to work out some of these on the graph. Partial products correspond to working through subgraphs.
Perhaps readers might find cases where working from the middle outward can be useful.
For example it would be possible to go from the middle outward using the Example of Figure \ref{hw}: we would go from  $c$ to $cd$ then
compute $acd$ and $bcd$.

\subsubsection{Edge Elimination}

\begin{figure}[H]
  \begin{center}
    \includegraphics[width=0.75\textwidth]{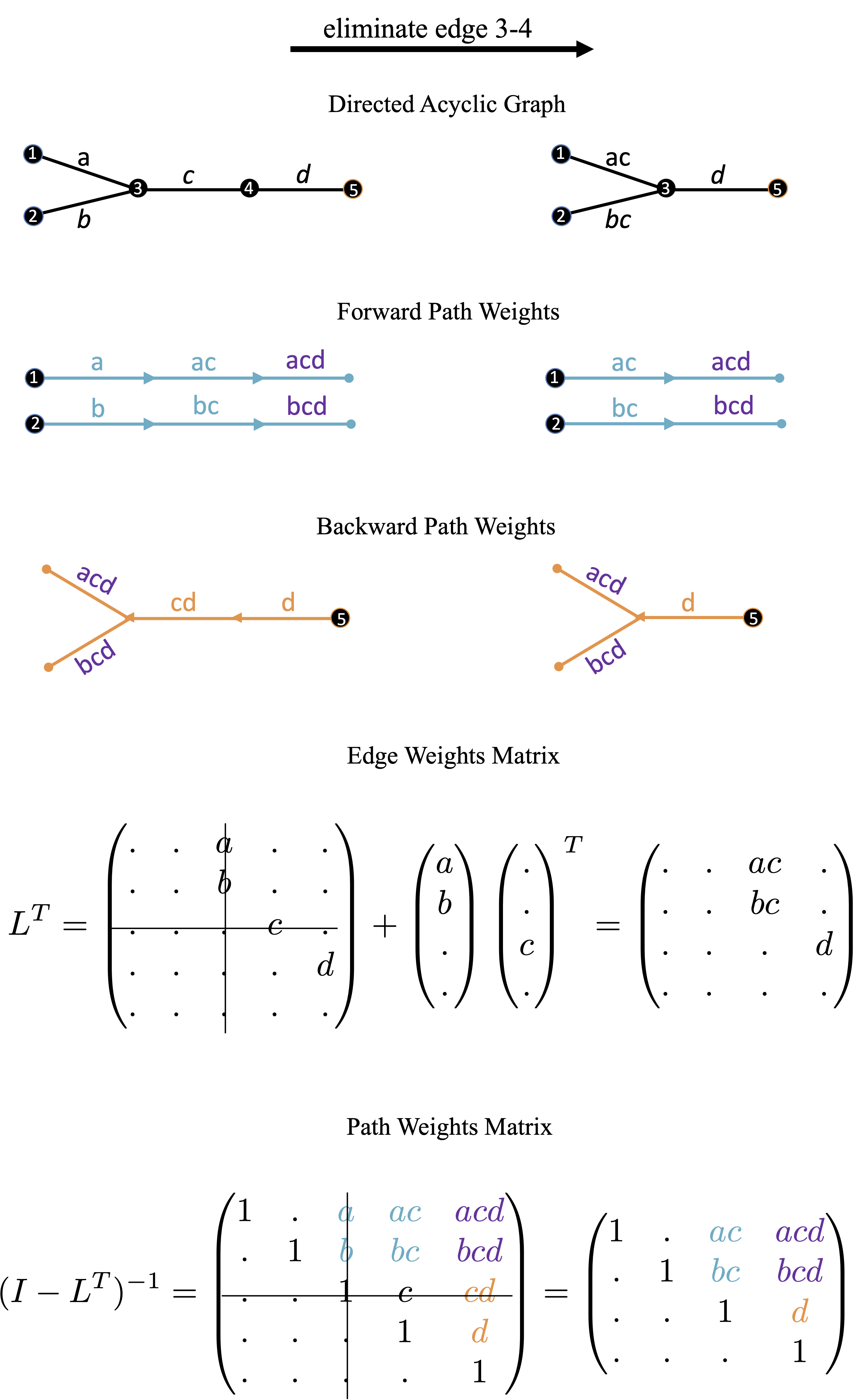}
  \end{center}
  \caption{\label{elim} Elimination of the edge from node 3 to node 4
    on the graphs and with matrices.  The matrix versions involve
    a rank one update to a row and column
    deleted matrix in the case of the Edge Weights Matrix and only a
    deletion of a row and column in the Path Weights Matrix.
  }
\end{figure}

It is possible to eliminate an edge
(and preserve the path weights)
by moving the weight of the edge to the
weights of the incoming edges.   We illustrate this in Figure \ref{elim}
by eliminating the edge from node 3 to node 4, moving the weight $c$ to the incoming edges by multiplication.
The corresponding linear algebra operation on $L^T$ is the deletion of column 3 and row 3
and the rank 1 update based on this column and row with the (3,3) element deleted.
The corresponding linear algebra operation on $(I-L^T)^{-1}$ is merely the deletion of column 3 and row 3.
This example  is representative of the general case.

\subsubsection{Edge addition at the Sink Node}

We will be  interested in the case where the edge weight graph is modified by adding one edge to the sink node.
Continuing our example from Figure \ref{elim}, we will  add
an edge ``$e$'' by starting with:
\begin{equation}
  \label{nowloss}
  \mbox{path weights} =
  \underbrace{\begin{pmatrix} 1 & 0  \\  0 & 1 \\ 0 & 0  \\  0 & 0 \\ 0 & 0 \end{pmatrix}^{\! T}}_{\mbox{sources}}
  \underbrace{  (I-L^T)^{-1}}_{   \tiny  \mbox{ path weights}  \mbox{ matrix} }
  \underbrace{\begin{pmatrix} 0 \\  0 \\ 0 \\  0 \\ 1 \end{pmatrix}}_{\mbox{sink}}
\end{equation}
and then updating by augmenting the graph with one end node to become
\begin{equation}
  \label{nowloss2}
  \mbox{updated path weights } =
  \underbrace{\begin{pmatrix} 1 & 0  \\  0 & 1 \\ 0 & 0  \\  0 & 0 \\ 0 & 0  \\ 0 & 0 \end{pmatrix}^{\! T}}_{\mbox{sources}}
  \underbrace{
    \begin{pmatrix}
      (I-L^T)^{-1} & . \\
      .            & 1
    \end{pmatrix}
    \begin{pmatrix}
      1 & . & . & . & . & .              \\
      . & 1 & . & . & . & .              \\
      . & . & 1 & . & . & .              \\
      . & . & . & 1 & . & .              \\
      . & . & . & . & 1 & \! \! e  \! \! \\
      . & . & . & . & . & 1              \\
    \end{pmatrix}}_{\small \mbox {updated} \mbox{ path weights}  \mbox{ matrix}}
  \underbrace{\begin{pmatrix} 0 \\  0 \\ 0 \\  0 \\ 0 \\ 1 \end{pmatrix}}_{\mbox{sink}} .
\end{equation}

The update from the path weights matrix in \eqref{nowloss} to the updated
path weights matrix in \eqref{nowloss2}
can  be verified in many ways.
One simple way is to look at the explicit elements of the path weights matrix
before and after and then notice that the new matrix has a column with one more element
$e$ augmented with a $1$.

It is an easy exercise in linear algebra to show that \eqref{nowloss2}  is the same as \eqref{newloss3} which
folds the added edge $e$  multiplicatively into the sink vector.
\begin{equation}
  \label{newloss3}
  \mbox{updated path weights } =
  \underbrace{\begin{pmatrix} 1 & 0  \\  0 & 1 \\ 0 & 0  \\  0 & 0 \\ 0 & 0 \end{pmatrix}^{\! T}}_{\mbox{sources}}
  (I-L^T)^{-1}
  \underbrace{\begin{pmatrix} . \\  . \\ . \\  . \\ e \end{pmatrix}}_{\mbox{sink}} .
\end{equation}

\subsection{Examples of DAGs and Weighted Paths}

\subsubsection{The ``Complete DAG'' and Weighted Paths}

\begin{figure}[H]
  \begin{center}
    \includegraphics[width=.9\textwidth]{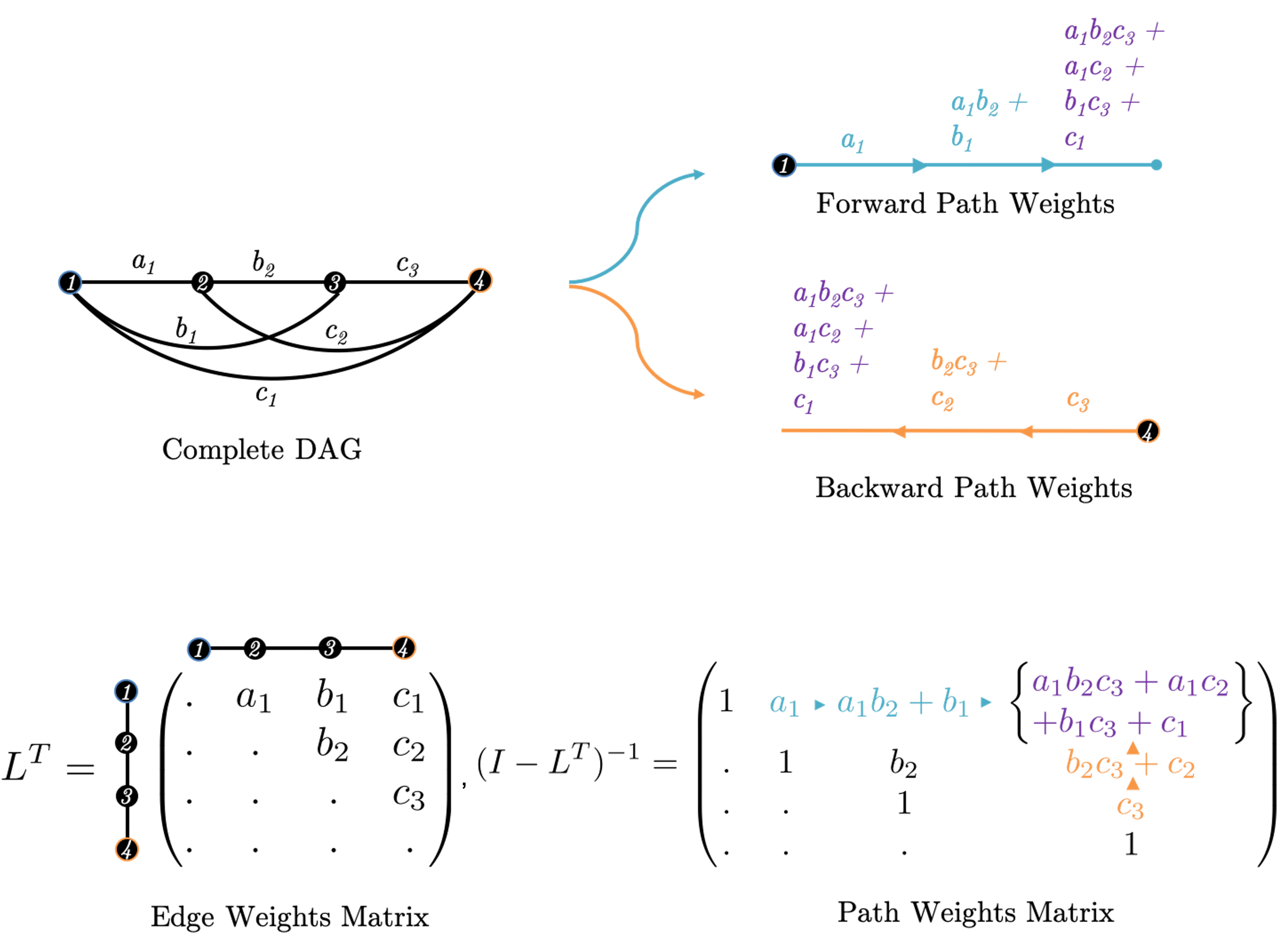}
    \caption{\label{fulldag}
      The complete DAG on four nodes illustrates a symmetric situation where forward and backward have the same complexity but arrive at the same answer through different operations.
    }
  \end{center}
\end{figure}

Consider as an example in Figure~\ref{fulldag}, the complete DAG on four nodes with graph weights evaluated through a forward and backward method. There is one source and one sink.  We find that this complete DAG example reveals most clearly the equivalence between path weights and the inverse matrix.

We see that  the forward path weights folds in the edges labelled ``$a$,'' then ``$b$,'' then ``$c$.''  This works through the matrix $L^T$ by columns. The backward mode folds in the edges with subscript ``3,'' then ``2,'' then ``1.''  This works through the matrix $L^T$ by rows from bottom to top.

\subsubsection{The ``multi-layer perceptron DAG''  and Weighted Paths}
\label{sec:mlp}

\begin{figure}[h]
  \begin{center}
    \includegraphics[width=\textwidth]{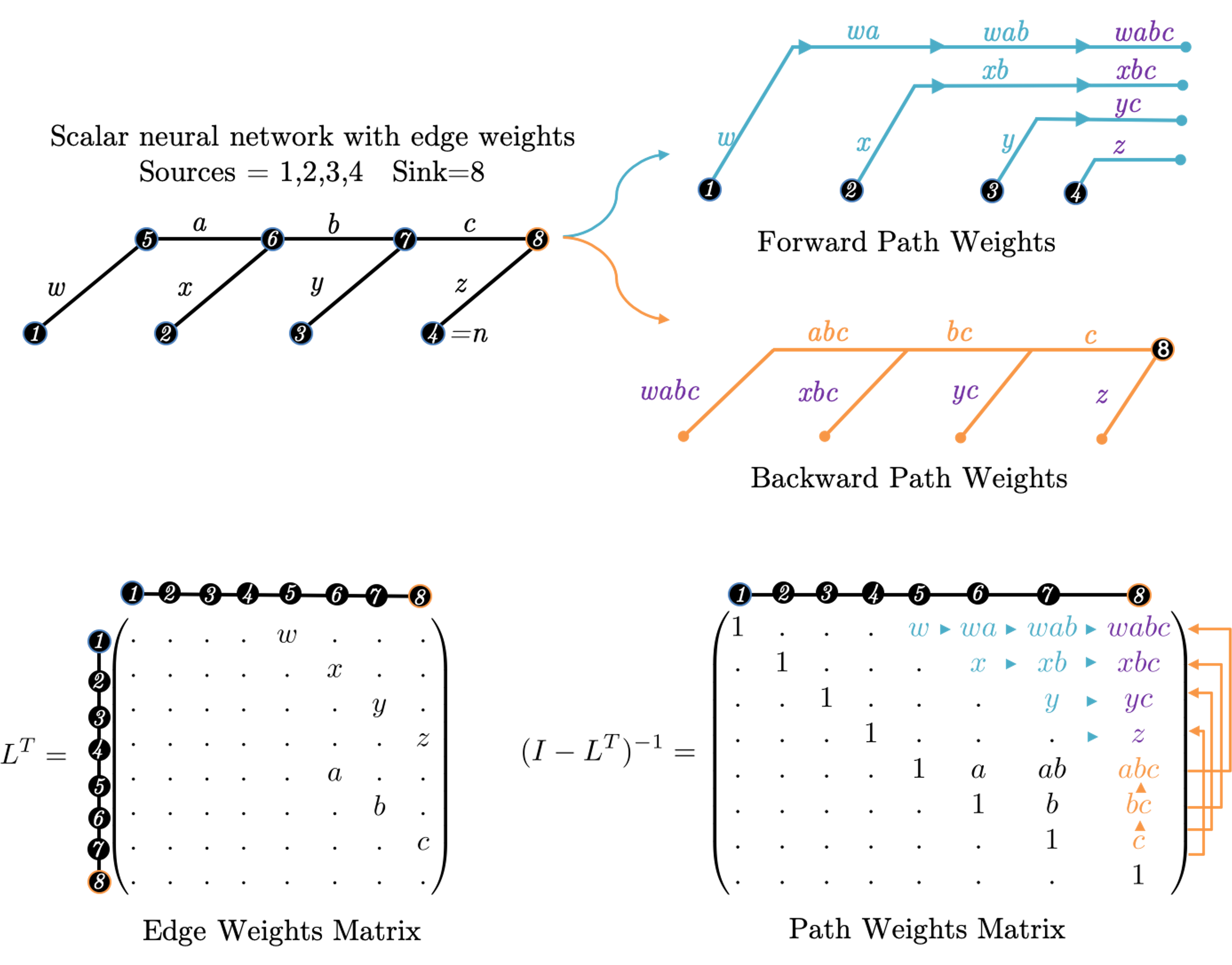}
    \caption{ This diagram contains most of what is needed to understand forward and
      backward propagation of derivatives through a MLP.
      The details of what the weights look like will come later.
      If we take
      $n=4$
      for the pictured network, the sources are labeled $1:n$ and the sink is labeled
      $2n$.
      Forward mode requires
      $n(n-1)/2$
      multiplications while backward mode requires
      $2n-3$
      multiplications.
    }
    \label{mlp}
  \end{center}
\end{figure}

Figure~\ref{mlp} is the DAG for the derivatives in a multi-layer perceptron (MLP). It may be thought of as a spine with feeds for parameters (nodes 1,2,3, and 4 in the figure).

If sources are labeled $1,\ldots,s$
(in Figure \ref{mlp}, $s=4$),
then  the top left $s$ by $s$ matrix in $L^T$ is the zero matrix
as there are no connections.
We can then write
\begin{equation}
  L^T  = \begin{pmatrix} 0 & M^T \\ 0 & \tilde{L}^T \end{pmatrix},
  \label{eqn:lt}
\end{equation}
where

\[
  M^T = \begin{pmatrix}
    w & . & . & . \\
    . & x & . & . \\
    . & . & y & . \\
    . & . & . & z
  \end{pmatrix},  
\tilde{L}^T = \begin{pmatrix}
    . & a & . & . \\
    . & . & b & . \\
    . & . & . & c \\
    . & . & . & . \\
  \end{pmatrix},
\]
where the matrix $M^T$ corresponds to connections between the sources and internal nodes, and $\tilde{L}^T$ corresponds to internal connections.  In this example $M^T$ is diagonal corresponding to a bipartite matching between nodes $1,2,3,4$ and $5,6,7,8$. The $\tilde{L}^T$ matrix represents internal connections, in this case it is the ``spine'' linearly connecting nodes $5,6,7,8$. 

Now we have
\[
  (I-L^T)  = \begin{pmatrix} I & -M^T \\ 0 & I-\tilde{L}^T \end{pmatrix},
  \mbox{ and } (I-L^T)^{-1} = \begin{pmatrix}  I &  M^T(I-\tilde{L}^T)^{-1}\\ 0 & (I-\tilde{L}^T)^{-1} \end{pmatrix}.
\]
If the last node is  the one  unique sink, then we obtain the useful formula
\begin{equation}
  \label{likeadjoint3}
  \mbox{Path weights} =  M^T(I-\tilde{L}^T)^{-1} \begin{pmatrix} 0 \\ \vdots \\0 \\ 1 \end{pmatrix}.
\end{equation}

We can now take a close look at Figure \ref{mlp} and fully grasp the path weight structure.
The spine consisting of $a,b,c$ and $1$ (understood) requires the computation of the cumulative suffix product  $1,c,bc,abc$. What follows is an element-wise multiplication by $z,y,x,w$,
from which we can calculate the last column of $M^T(I-\tilde{L}^T)^{-1}$
\begin{equation}
  \label{likeadjoint}
  M^T(I-\tilde{L}^T)^{-1} \begin{pmatrix} 0 \\ 0 \\ 0 \\ 1 \end{pmatrix}=
  \begin{pmatrix} wabc \\ xbc \\ yc \\ z \end{pmatrix}.
\end{equation}

\subsection{Computational Graphs, Derivative Graphs, and their superposition}

Many treatments of automatic differentiation introduce
computational graphs at the start of the
discussion.  Our treatment shows that this is not necessary. However, in the end the key application of edge weights will be as derivatives of computed quantities.  To this end, we define

\begin{defn}
  A computational graph is a node labelled DAG, where leaf nodes consist of variable names, and non-leaf nodes contain variable names and
  formulas that depend on incoming variables.
\end{defn}

We remark that there are variations on where the variable names and formulas live on a computational graph, but we believe
the definition here is the cleanest when wishing to incorporate derivative information.

\begin{figure}[h]
  \begin{center}
    \includegraphics[width=\textwidth]{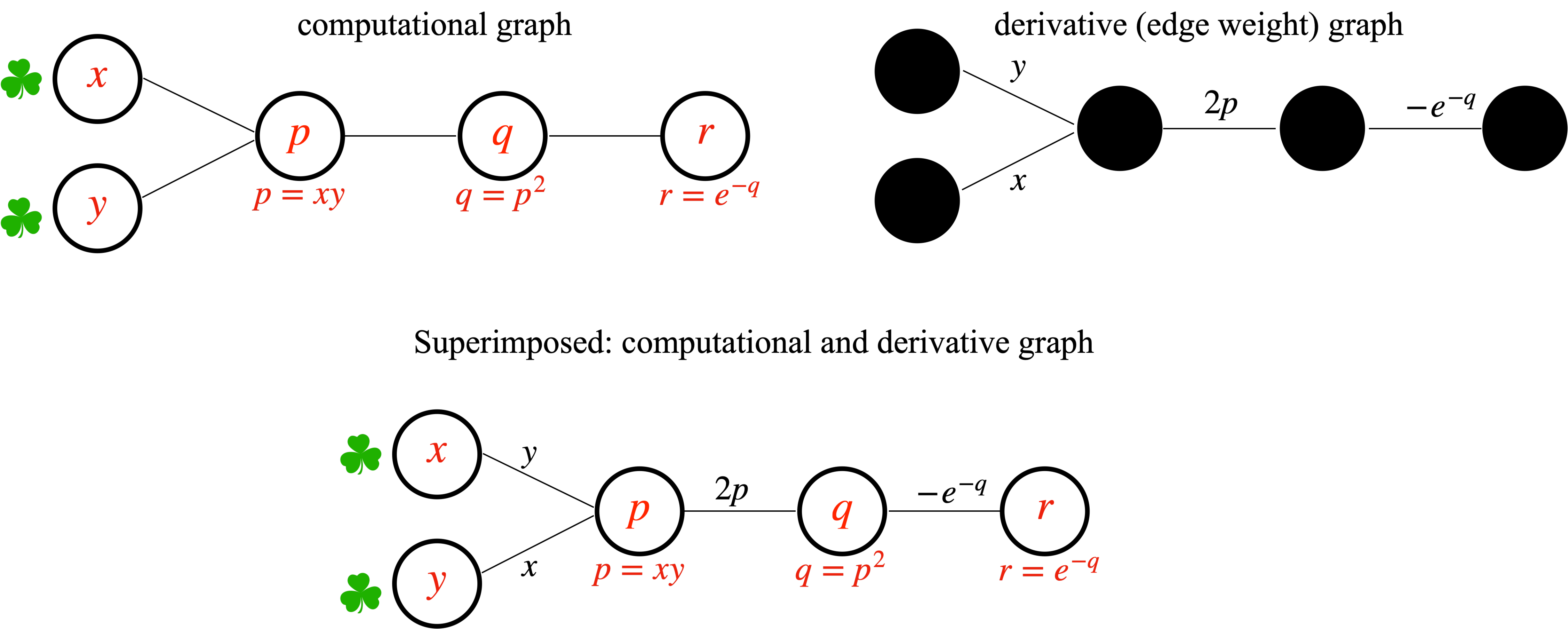}
    \caption{An example of Computational (node) graph, derivative (edge) graph, and their ``superposition.''
    }
    \label{figcomp}
  \end{center}
\end{figure}

\begin{table}
\begin{center}
  \begin{tabular}{ m{3cm}| m{8cm} }
    \toprule
    Player                            & Description                                          \\
    \midrule
    Edge weight from node $i$ and $j$ & These are the derivatives of  one step of
    a computation.  These can be scalars
    but in general these are Jacobian matrices (or operators).                               \\
    \midrule
    Path weight from node $i$ to  $j$ & These are the derivatives that
    reach back into a chain of computations.
    The chain rule states that if you multiply
    (``chain together'') the  derivatives at each step
    you get the  dependence of one variable on an earlier variable.                          \\
    \midrule
    Source                            & The sources in the graph are typically
    parameters in real computations, as many modern applications
    are interested in the derivatives with respect to the input
    parameters.                                                                              \\
    \midrule
    Sink                              & The sink is usually what is known as a loss function
    in modern applications.                                                                  \\
    \bottomrule
  \end{tabular}
  \caption{A dictionary translating graph elements to AD (cf. Figure~\ref{figcomp}).}
\end{center}
\end{table}

\subsubsection{The chain rule, derivatives, and Jacobians}

Here we say explicitly how the edge weights and path weights
relate to derivatives in a computation.

Consider the computation from Figure \ref{figcomp}, the next three
algorithms show the computation, the derivatives of each
line of code, and the overall derivatives.  We see that
the one step derivatives are edge weights and the overall derivatives
are path weights.

If the final output is a scalar, we immediately have that the gradient with respect to the source $x$ and $y$ (in Figure \ref{figcomp}) is exactly the path weight defined in \eqref{likeadjoint3},
\begin{equation}
  \label{likeadjoint2}
  \text{gradient = the last column of }\ M^T(I-\tilde{L}^T)^{-1}, 
\end{equation}
which corresponds to the output in Algorithm 3.8 with
\[
M^T = \begin{pmatrix}
  y & \cdot & \cdot\\
  x & \cdot & \cdot\\
\end{pmatrix}, \quad 
\tilde{L} = \begin{pmatrix}
  \cdot & 2p & \cdot\\
  \cdot & \cdot & -e^{-q}\\
  \cdot & \cdot & \cdot\\
\end{pmatrix}.
\]
 Equation~\eqref{likeadjoint2} fully describes backpropogation. For completeness,  the term "forward propagation" (forward pass) describes the process of executing the computational graph in a forward direction (left to right), storing the intermediate values that are subsequently utilized in backpropagation (backward pass).

\begin{minipage}{0.46\textwidth}
  \begin{algorithm}[H]
    \caption{Simple Algorithm Example from Figure \ref{figcomp}}\label{gen1}
    \begin{algorithmic}[1]
      \State $ p \gets \text{multiply}(x,y)$
      \State $q \gets \text{square}(p)$
      \State $r \gets  \text{exp\_neg}(q)$
      \State {\color{ekinorange} output $r$}  \end{algorithmic}
  \end{algorithm}
\end{minipage}
\hfill
\begin{minipage}{0.46\textwidth}
  \begin{algorithm}[H]
    \caption{Edge weights
      (derivatives of one line of code)}
    \label{gen2}
    \begin{algorithmic}[1]
      \State $d\{\text{multiply}(x,y)\}/d x = y$ \ \ \ ($=\frac{\partial p}{\partial x} $)
      \State $d\{\text{multiply}(x,y)\}/d y  = x$ \ \ \ ($=\frac{\partial p}{\partial y} $)
      \State $d\{\text{square}(p)\}/d p = 2p$ \ \ \  \ \ \ \ \ ($=\frac{\partial q}{\partial p} $)
      \State $d\{\text{exp\_neg}(q)\}/d r =-e^{-q}$ \ \ ($=\frac{\partial r}{\partial q} $)
    \end{algorithmic}
  \end{algorithm}
\end{minipage}

\begin{center}
  \begin{minipage}{.9\textwidth}
    \begin{center}
      \begin{algorithm}[H]
        \caption{Path weights
          (Chained derivatives)
        }\label{gen3}
        \begin{algorithmic}[1]
          \State $d r/d x = y \times 2p \times (-e^{-q})$
          (Chain  lines 1,3, and 4 of Algorithm \ref{gen2})
          \State $d r/d y = x \times 2p \times (-e^{-q})$
          (Chain lines 2,3, and 4 of Algorithm \ref{gen2} )
        \end{algorithmic}
      \end{algorithm}
    \end{center}
  \end{minipage}
\end{center}

\section{Linear Operators as elements of Matrices}
\label{sec:matrix}

\label{bigops}

We will illustrate in Section \ref{challenge} the value of software that allows
linear operators as elements of matrices.  Here we set
the mathematical stage, starting with a question.\\

\emph{Consider a matrix transformation of $X$ such as
$T_{A,B}: X \mapsto BXA^T,$ how should we represent the Jacobian $\partial T_{A,B}/\partial X$?}\\

\noindent Before we answer, we remind the reader how the Kronecker product works.  One view of the Kronecker product
$A \otimes B$ of two matrices is that it multiplies every element in $A$ times every element of $B$ placing the elements in such
a way that we have the identity
\[ (A \otimes B) \text{vec}(X) = \text{vec}(BXA^T),\]
where  $vec$ denotes the flattening of a matrix $X$ into a vector by stacking its  columns. We may abuse notation when there is no confusion and write
\[ (A \otimes B) (X) = BXA^T,\]
for the linear operator $T_{A,B}$ that sends $X$ to $BXA^T$. Identifying the matrix $A \otimes B$ with the operator is more than a handy convenience, it makes computations practical in a software
language that allows for this. Table~\ref{tab:operators2} defines some operators of interest.

\begin{table}[h]
  \centering
  \begin{tabular}{m{3cm}|c|l|lc}

    \toprule
                                  & Symbol          & Definition                         & \multicolumn{2}{c}{ Dense Representation}                      \\
    \midrule
    Kronecker Product of $A,B$    & $A \otimes B$   & $X \mapsto  BXA^T\rule{0pt}{3ex} $ & $A \otimes B$                             & $m_1n_1 \times mn$ \\
    \midrule
    Left Multiplication by $B$    & $B_L$           & $ X \mapsto BX$                    & $I \otimes B$                             & $m_1n \times mn$   \\
    \midrule
    Right Multiplication by $A$   & $A_R $          & $X \mapsto XA$                     & $A^T \otimes I$                           & $m n_1 \times mn$  \\
    \midrule
    Hadamard Product with $M$     & $M_H$           & $X \mapsto M .*X $                 & diag(vec($M$))                            & $m n \times mn$    \\
    \midrule
    Matrix inner product with $G$ & $G^{T_\bullet}$ & $X \mapsto \mbox{Tr}(G^TX) $       &
    vec($G$)$^T$                  & $1 \times mn$                                                                                                         \\
    \midrule
  \end{tabular}
  \caption{ Matrix Operators  and the size of their  dense representations
    assuming $X: m \times n,$ \ \ $A:n_1 \times n,$  \  \   $B:m_1 \times m,$ \ \  $M:m \times n,$ and $G:m \times n$. We overload $A \otimes B$ to be both the operator and the matrix.
  }
  \label{tab:operators2}
\end{table}

Consider the inner product (matrix dot product) $\langle X,Y \rangle =\mbox{Tr}(X^TY)$.
The identity $\langle X,AY \rangle=\langle A^TX,Y \rangle$
implies $(A_L)^T=(A^T)_L,$
in words, the  operator adjoint with respect to the operator $A_L$ (left multiplication by $A$) is left multiplication by $A^T$. The  operator transposes are $(A_L)^T = (A^T)_L$, $(B_R)^T = (B^T)_R$, and $(M_H)^T = M_H$ (symmetric).

We wish to propose a carefully
thought out notation for another useful operator,  $G^{T_\bullet}$ (``$G$ transpose dot''), the matrix inner (or dot) product with $G$.
\begin{defn}
  Let $G^{T_\bullet}$ (``$G$ transpose dot'') denote the matrix inner (or dot) product with $G$. This operator takes a matrix $X$ of the same size as $G$ and returns the scalar, $G^{T_\bullet}X :=$ \text{Tr}$(G^TX)$= vec($G$)$^T$vec(X) = $\sum_{i,j} G_{ij}X_{ij}$.
\end{defn}

Many communities choose a notation where small Roman letters
denote a column vector, so that $x \mapsto g^Tx$ denotes
a linear function of $x$.  Those who are used to this notation
no longer ``see'' the transpose so much as turning a column into
a row, but rather  they see the linear function $g^T$
as an object that acts on (``eats'') vectors and returns scalars.
In the same way we propose that one might denote a linear
function of a matrix $X\mapsto \mbox{Tr}(G^TX)$ with the operator
notation $X \mapsto G^{T_\bullet}X$, an operator that ``eats''
matrices and returns scalars.

\begin{lemma}
  If  the superscript ``$()^T$''  is overloaded to denote real operator adjoint or matrix transpose
  as appropriate, ${\cal L}$ is a linear
  operator and $G$ is a matrix, then we have the operator identity:
  $({\cal L}^TG)^{T_\bullet}=G^{T_\bullet}{\cal L}.$
  Notice that if we pretend all letters are just matrices and
  you ignore the dot, the notation has the appearance of the familiar transpose rule.

\end{lemma}
\begin{proof}
  We have that for all $X$,
  \[
    ({\cal L}^TG)^{T_\bullet}X = \langle {\cal L}^TG,X \rangle  = \langle G, {\cal L}X \rangle =
    G^{T_\bullet} {\cal L}X,
  \]
  showing that as operators
  $({\cal L}^TG)^{T_\bullet}=
    G^{T_\bullet} {\cal L}.$ 
  \end{proof} As an example, we have
  \[(A^T_LG)^{T_\bullet} = X \mapsto \mbox{Tr}( (A^TG)^TX), \quad \text{and}\quad G^{T_{\bullet}}A_L = X \mapsto  \mbox{Tr}(G^TAX),\]
  which shows that $(A^T_LG)^{T_\bullet} = G^{T_{\bullet}}A_L$. We encourage the reader to follow the matrices
  $A,G,A^T$ and the operators
  $A_L^T,A_L,(A^T_LG)^{T_\bullet},G^{T_{\bullet}}.$
  (See Section \ref{notation} for why this notation can be valuable.)

\section{Operator Methodology}

We proceed from matrices of scalars to
matrices of vectors to matrices of operators in Sections 5.1, 5.2, and 5.3.
ultimately taking advantage of Julia's
capabilities.
We encourage the reader to compare
the matrices in each of these sections.
Section 5.4 illustrates the power of
the $G^{T_\bullet}$ notation, while
Section 5.5 shows the relationship to the
adjoint method that is well known in
the field of scientific computing.

\begin{table}
  \centering
  \resizebox{\textwidth}{!}{\begin{tabular}{m{1.8cm}|l llll}
    \toprule
                                                      & \multicolumn{5}{c}{\includegraphics[width=.4\textwidth]{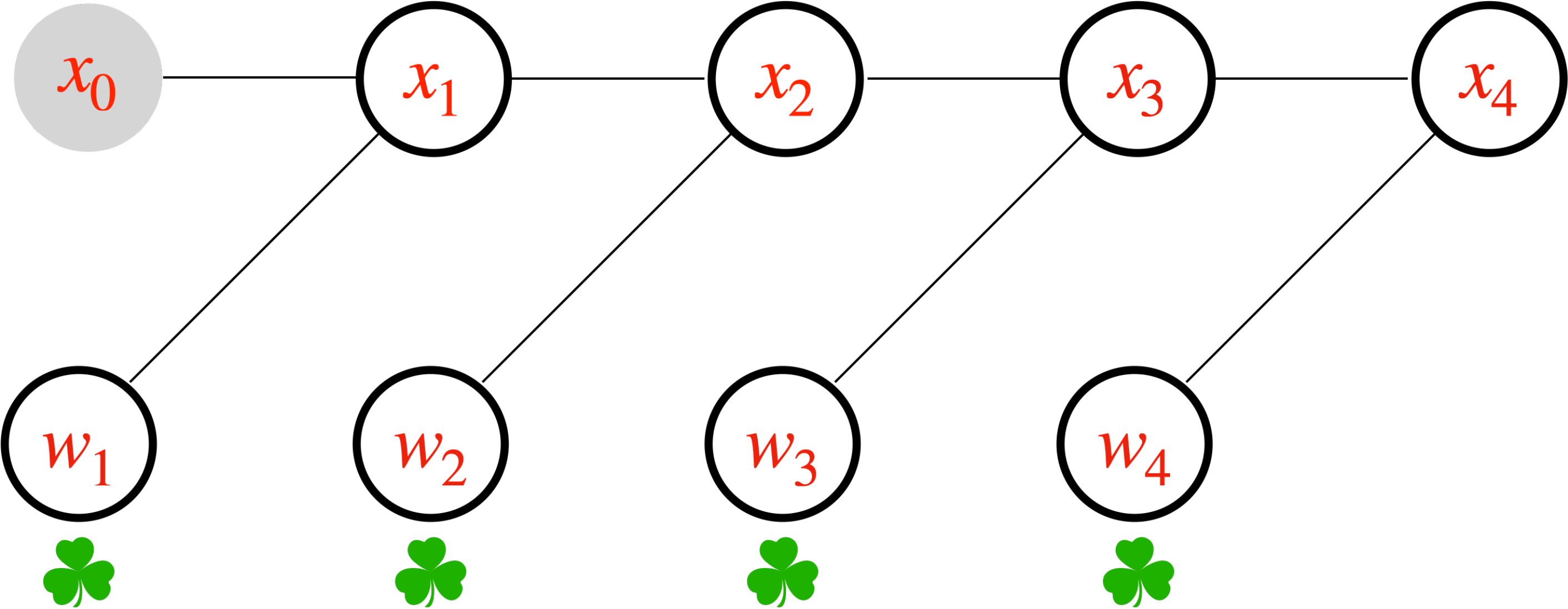}}                                                                                                                                                                                                                                                                                        \\
    \midrule
    Gradient w.r.t. parameters $p$ (leaf nodes)             & $\nabla_p\mathcal{L}= $  \hspace{.2in} \textcolor{red}{$M^T$}                       & $\times$                                                                                                              & \hspace{.3in}   \textcolor{blue}{$(I - L)^{-1}$}                                                  & $\times$ & \hspace{.1in}  \textcolor{cyan}{$g$} \\ \\
                                                            & = $\color{red}{\begin{pmatrix}
                                                                                   m_1 &     &     &     \\
                                                                                       & m_2 &     &     \\
                                                                                       &     & m_3 &     \\
                                                                                       &     &     & m_4
                                                                                 \end{pmatrix}}^T$                                                      & $\times$                                                                                                              &
    $\color{blue}{\begin{pmatrix}
                        I    &      &      &   \\
                        -l_2 & I    &      &   \\
                             & -l_3 & I    &   \\
                             &      & -l_4 & I \\
                      \end{pmatrix}}^{-T}$
                                                            & $\times$                                                                            & $\color{cyan}{\begin{pmatrix} . \\ . \\ . \\ g_4 \end{pmatrix} }$                                                                                                                                                                                                           \\
    \midrule
    (i) Scalar                                              &                                                                                     &                                                                                                                       &                                                                                                   &                                                 \\
    $p = \{w_i\}$                                           & \multicolumn{2}{l}{$m_i = \delta_ix_{i-1}$}                                         & \multicolumn{2}{l}{$l_i = \delta_{i}w_i$}                                                                             & $g_4 = {\mathcal L}^\prime(x_4)$                                                                                                                    \\ \\
    (ii) Vector                                             &                                                                                     &                                                                                                                       &                                                                                                   &                                                 \\
    $p = \{[w_i, b_i]\}$                                    & \multicolumn{2}{l}{$m_i = [\delta_ix_{i-1} \  \delta_i]$}                           & \multicolumn{2}{l}{\hspace*{.1in} \textquotesingle \textquotesingle \hspace*{.1in} \textquotesingle \textquotesingle} & \hspace*{.1in} \textquotesingle \textquotesingle \hspace*{.1in} \textquotesingle \textquotesingle                                                   \\ \\
    (iii) Matrices                                          &                                                                                     &                                                                                                                       &                                                                                                   &                                                 \\
    $p=\{[W_i, B_i]\}$                                      & \multicolumn{2}{l}{$m_i =
    [{\Delta_i}_H  \circ {X_{i-1}}_R \,\, \ {\Delta_i}_H]$} & \multicolumn{2}{l}{$l_i =  {\Delta_{i}}_H \circ {W_{i}}_L$}                         & $g_4 = \nabla_{X_4} {\cal L}$                                                                                                                                                                                                                                               \\
                                                            & \multicolumn{4}{c}{${}^\uparrow \! \! \! - {\rm Operators}{} - \! \! \! ^\uparrow$} &                                                                                                                                                                                                                                                                             \\
    \bottomrule
  \end{tabular}}
  \caption{Algebraic Structure for an MLP when the parameters (i.e. the set of leaf nodes collectively referred to as $p$) are (i) only scalar weights (ii) a weight/bias vector, and (iii) a vector of weight/bias matrices.  We emphasize the common algebraic structure and the benefit of software that can represent matrices of vectors and matrices of operators.}
  \label{tab:my_label}
\end{table}

\subsection{Matrices of scalars}

\begin{algorithm}[H]
  \caption{Scalar MLP without Bias (forward propagation)}\label{nnalg}
  \begin{algorithmic}[1]
    \State {\color{ekinorange} Input data $x_0$, initial weights $w_i, \quad i = 1, \cdots, N$}
    \State Select activation functions $h_i(\cdot)$ such as sigmoid, tanh, ReLU, etc.
    \For{$i$ = 1 to $N$}
    \State $x_{i}  \gets h_i(w_{i}  x_{i-1} $)
    \State {\color{ekinblue}$(\delta_i \gets h_i^\prime(w_{i}  x_{i-1}$)) }
    \EndFor
    \State {\color{ekinorange} Output $x_{N}$}
  \end{algorithmic}
\end{algorithm}

The simple case of scalar neural networks (shown in Algorithm~\ref{nnalg}) without bias shows the power of
the graph approach.  However, the full power is revealed
in the coming sections.  Here we remind the reader
of the algorithm, draw the graphs, and instantly
write down the linear algebra that provides the gradients
through backpropogation. (The graphs and matrices
are illustrated for $N=4$ for ease of presentation.) 
\begin{figure}[H]
  \begin{center}
    \includegraphics[width=\textwidth]{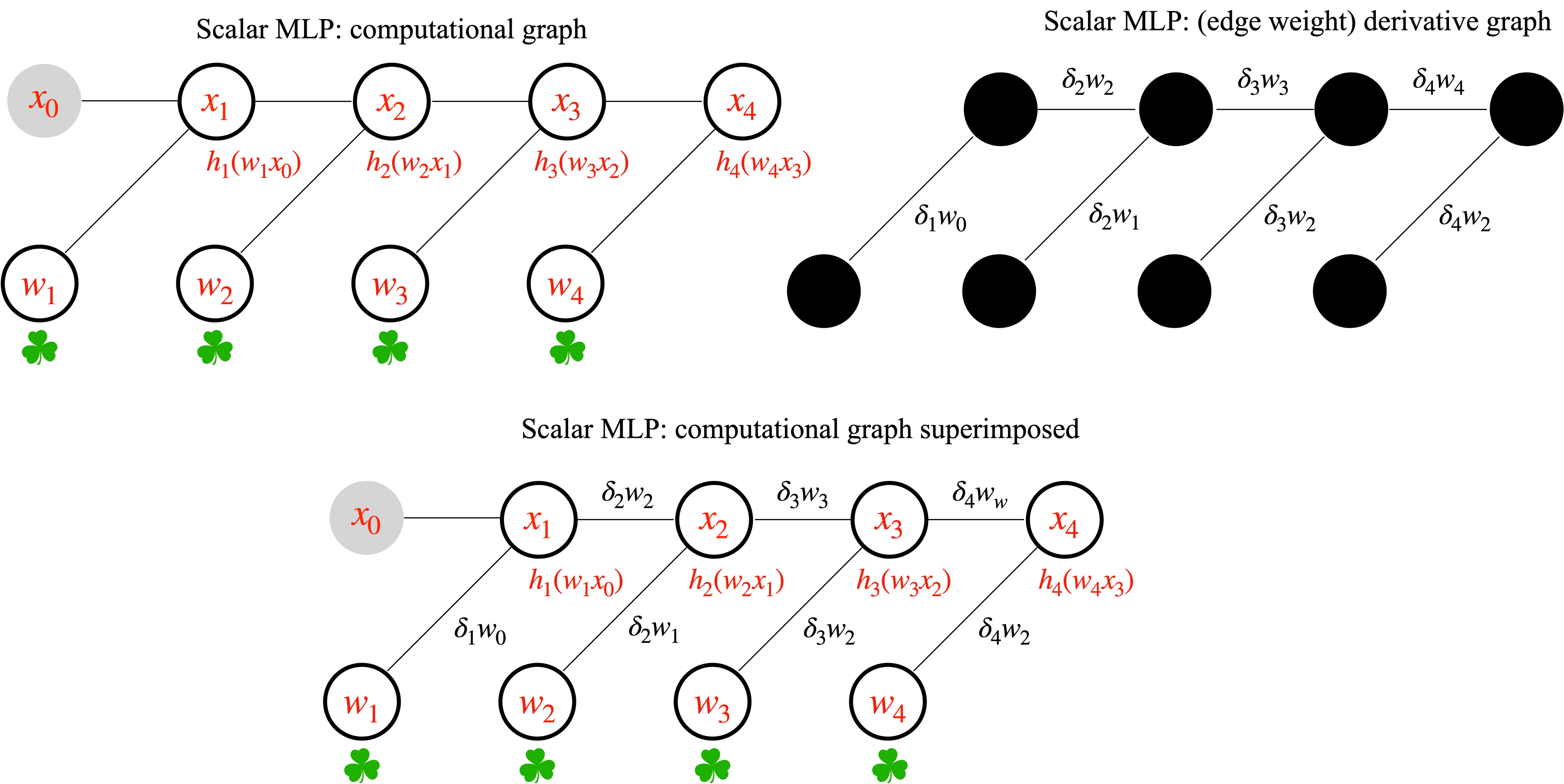}
    \caption{Top left: computational graph of a scalar MLP. This computation, which has nothing to do with derivatives, is often referred to as forward propagation because of its direction. Evaluation must generally necessarily go from left to right. Top right: derivative edge weights. Since derivatives are linear, multiple directions are possible to evaluate the products. Bottom: the superimposed graph showing both the forward computation and the derivative edge weights.}
    \label{dnn}
  \end{center}
\end{figure}
Starting with
\[\tilde{L}^T = \begin{pmatrix}
    . & \delta_2 w_1 & .            & .            \\
    . & .            & \delta_3 w_2 & .            \\
    . & .            & .            & \delta_4 w_3 \\
    . & .            & .            & .            \\
  \end{pmatrix},
  M = \begin{pmatrix}
    \delta_1 x_0 & .            & .            & .             \\
    .            & \delta_2 x_1 & .            & .             \\
    .            & .            & \delta_3 x_2 & .             \\
    .            & .            & .            & \delta_4  x_3
  \end{pmatrix},
\]
it is an immediate consequence of our graph theory
methodology which concluded with
\eqref{newloss3}
and
\eqref{likeadjoint}
that the backpropagated gradient is computed by evaluating efficiently

$$\nabla_w {\cal L} =
  \begin{pmatrix}
    \delta_1 x_0 &              &              &               \\
                 & \delta_2 x_1 &              &               \\
                 &              & \delta_3 x_2 &               \\
                 &              &              & \delta_4  x_3
  \end{pmatrix}
  \begin{pmatrix}
    1             &               &               &   \\
    -\delta_2 w_1 & 1             &               &   \\
                  & -\delta_3 w_2 & 1             &   \\
                  &               & -\delta_4 w_3 & 1 \\
  \end{pmatrix}^{-T}
  \begin{pmatrix} . \\ . \\ . \\ {\mathcal L}^\prime(x_4) \end{pmatrix}
$$

\subsection{Matrices of vectors}

As a baby step towards the matrices of operators approach,
we show how one can (optionally) group weights and biases
that appear in a neuron.  Algorithm \ref{nnalg} is modified so that
$w_i x_{i-1}$ is replaced with $w_i x_{i-1}+b_i$.
In the interest of space, we
will simply write the answer of $\nabla_{[w,b]} {\cal L}$ and discuss its format,
$$
  \begin{pmatrix}
    [{\delta_1}  {x_0} \hspace*{.05in} \delta_1]                                                                                                                                                                                          \\
     & \hspace{-.2in}[{\delta_2}  {x_1}  \hspace*{.05in}  \delta_2]                                                                                                                                                                       \\
     &                                                              & \hspace{-.2in}\hspace{.05in} \  [{\delta_{3}}  {x_{2}} \hspace*{.05in}  {\delta_{3}}]                                                                               \\
     &                                                              &                                                                                       & \hspace{-.2in} \  [{\delta_{4}}  {x_{3}}  \hspace*{.05in}  {\! \delta_{4}}]
  \end{pmatrix}^T \! \! \!
  \begin{pmatrix}
    1             &               &               &   \\
    -\delta_2 w_1 & 1             &               &   \\
                  & -\delta_3 w_2 & 1             &   \\
                  &               & -\delta_4 w_3 & 1 \\
  \end{pmatrix}^{\! -T} \! \! \!
  \begin{pmatrix} .\\ . \\ . \\ {\mathcal L}^\prime(x_4) \end{pmatrix}  \! \! .
$$

We see we have an ordinary matrix back substitution followed by multiplication
by a diagonal matrix of  row vectors of length 2
so that the result is
a vector of column vectors of length 2 which nicely packages the gradients
with respect to the weight and bias in each neuron.
We remark that the transpose applies recursively
in the diagonal matrix.  The transpose is overkill
in this case but is critical in the next section.

\subsection{Matrices of operators}

Letting ${\cal I}$ denote the identity operator and
empty space the zero operator, we have the following
\begin{displaymath}
  \begin{split}
    \nabla_{[W,B]} {\cal L} &= \begin{pmatrix}
      [{\Delta_1}_H  \circ {X_0}_R \,\, {\Delta_1}_H]                                                                                                                                                                       \\
       & [{\Delta_2}_H  \circ {X_1}_R \,\, {\Delta_2}_H]                                                                                                                                                                    \\
       &                                                 &  & \hspace{-.2in} \  [{{\Delta_{3}}_H}  \circ {{X_{2}}_R} \,\, {{\Delta_{3}}_H}]                                                                                 \\
       &                                                 &  &                                                                               & \hspace{-.2in} \  [{{\Delta_{4}}_H}  \circ {{X_{3}}_R} \,\, {{\Delta_{4}}_H}]
    \end{pmatrix}^T\\
    &\qquad\qquad\quad\quad\ \ \times
    \begin{pmatrix}
      {\cal I}                    &                                     &                                                 \\
      -{\Delta_2}_H \circ {W_2}_L & {\cal I}                            &                                                 \\
                                  & -{{\Delta_{3}}_H} \circ {{W_{3}}_L} & {\cal I}                             &          \\
                                  &                                     & - {{\Delta_{4}}_H} \circ {{W_{4}}_L} & {\cal I} \\
    \end{pmatrix}^{-T}
    \begin{pmatrix}
      . \\. \\ . \\
      \nabla_{X_4} {\cal L}
    \end{pmatrix}
  \end{split}
\end{displaymath}
for the matrix neural network in Algorithm
\ref{matrixneuralnet}. The entries of our matrix of operators may be read immediately from
the differential of line 4 of Algorithm \ref{matrixneuralnet}:
\begin{align*}
  dX_{i} &= d\left[h_i(W_iX_{i-1} + B_i\right]\\
  &=   ({\Delta_i}_H  \circ {{}X_{i-1}}_R) dW_i \, +  {\Delta_i}_H dB_i
  + ({\Delta_i}_H \circ {W_i}_L )dX_{i-1},
\end{align*}
where $\Delta_i$ is the gradient matrix, and the definitions of the operators ${\Delta_i}_H$, ${W_i}_L$, and ${X_{i-1}}_R$ are given in Table~\ref{tab:operators2}.

\begin{algorithm}[H]
  \caption{Matrix MLP (forward propagation)}\label{matrixneuralnet}
  \begin{algorithmic}[1]
    \State {\color{ekinorange} Input data $X_0$ ($n_0 \times k$)}, and initial weight matrices and corresponding bias terms {\color{ekinorange} $W_i (n_i \times n_{i-1}), B_i (n_i \times k)$}
    \State Select activation functions $h_i(\cdot)$ such as sigmoid, tanh, ReLU, etc.
    \For{$i$ := 1 to $N$}
    \State $X_{i} \gets h_i(W_i*X_{i-1} + B_i)$.
    \State (${\color{ekinblue} \Delta_i \gets h_i^\prime(W_i*X_{i-1} + B_i)}$)
    \EndFor    
    \State {\color{ekinorange} output $X_{N}$}
  \end{algorithmic}
\end{algorithm}

\subsection{The Power of Notation}
\label{notation}

\begin{figure}[H]
  \begin{center}
    \includegraphics[width=\textwidth]{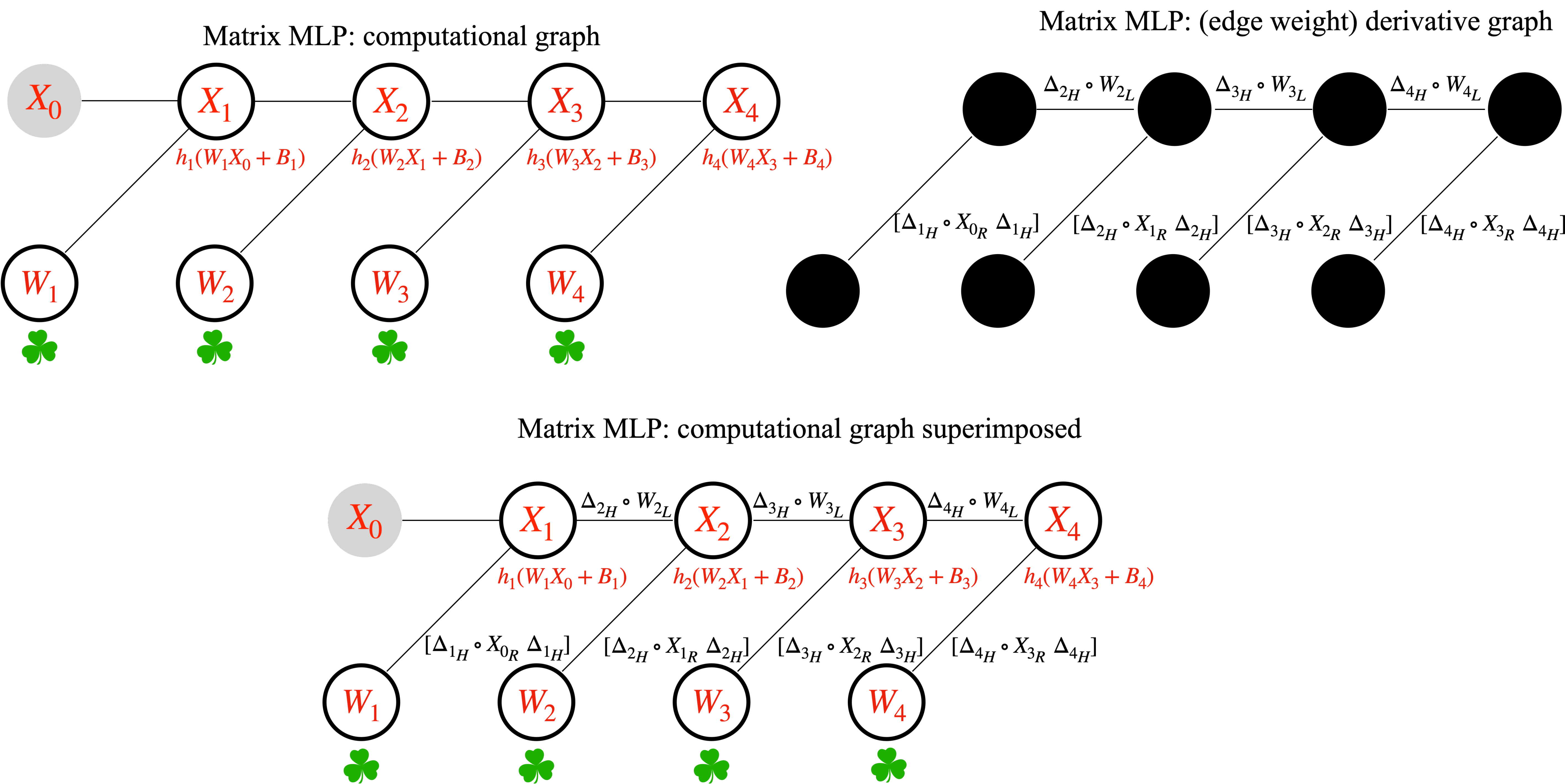}
    \caption{Computational and derivative graphs of a matrix MLP with their superimposed version. Compared to Figure~\ref{dnn} of the scalar MLP, everything remains the same except for two changes: the elements in the computational graph are now matrices, and the edges in the derivative graph have been replaced by operators.}
    \label{mat_dnn}
  \end{center}
\end{figure}

We read directly off the edge weight graph in Figure~\ref{mat_dnn} that for a matrix neural network we have
        \begin{equation}
          \textsc{Forward Mode Operators (right to left)}\nonumber\\          
          \begin{split}
          \frac{\partial \mathcal{L}}{\partial W_i} &=    G^{T_\bullet}
          (\Delta_N)_H(W_{N})_L\ldots(W_{i+2})_L(\Delta_{i+1})_H(W_{i+1})_L(\Delta_i)_H(X_{i-1})_R  \\ 
          \frac{\partial \mathcal{L}}{\partial B_i} &=  \
          G^{T_\bullet}
          (\Delta_N)_H(W_{N})_L\ldots(W_{i+2})_L(\Delta_{i+1})_H(W_{i+1})_L(\Delta_i)_H
          \label{fmo2}
          \end{split}
        \end{equation}
or going the other way we have,
        \begin{equation}
          \textsc{backward Mode Operators (right to left)}\nonumber\\   
          \begin{split}
          \left[\frac{\partial \mathcal{L}}{\partial W_i}\right]^{T_\bullet} &=
          \left\{ (X_{i-1}^T)_R (\Delta_i)_H (W_{i+1}^T)_L(\Delta_{i+1})_H(W_{i+2}^T)_L \ldots (W_{N}^T)_L(\Delta_N)_H G \right \}^{T_\bullet} \\
          \left[ \frac{\partial \mathcal{L}}{\partial B_i}\right]^{T_\bullet} &=  \
          \left\{  (\Delta_i)_H (W_{i+1}^T)_L(\Delta_{i+1})_H(W_{i+2}^T)_L \ldots (W_{N}^T)_L(\Delta_N)_H G \right \}^{T_\bullet}.
          \end{split}
        \end{equation}

\textbf{Understanding these operators.} The forward operators in Equation (\ref{fmo2})
may be thought of as sensitivity operators or as a means of computing
the full gradient.  As a sensitivity operator, one can
state that the directional derivative of $\mathcal{L}$
in the direction
$ \Delta W_i$ is $\frac{\partial \mathcal{L}}{\partial W_i}( \Delta W_i$).
  Alternatively, each operator can be written out as a (large) matrix,
  and ultimately a gradient can be computed. The backward operator is intended to be evaluated from right to left
  inside the braces.  Doing so computes the gradient directly.

  We hope the reader appreciates the power of the ``$T_\bullet$''
  notation, whereby one feels we are taking transposes of matrices
  and reversing order, but in fact we are transposing the operators.
  Either way the operators can be read right off the graphs.

\subsection{Relationship to the Adjoint Method of scientific computing}

  We will show how to derive \eqref{likeadjoint2} and \eqref{likeadjoint} using the adjoint method
  so-named because of its focus on the transpose (the adjoint) of the Jacobian. We encourage interested readers to  see
  \cite{johnson2007notes} and  \cite{bradley2010pde} to learn more about adjoint methods in numerical computation.

  We find it satisfying that the graph theoretic interpretation of backward mode AD and the adjoint method of scientific computing yield the same answer from two very different viewpoints.

  Consider a  general computation with  known constant  input $x_0 \in \mathbb{R}$ and parameters $p = [p_1,\ldots,p_k]$:

  \begin{center}
    \begin{algorithm}[H]
      \caption{General Computation}\label{gen}
      \begin{algorithmic}[1]
        \State {Input constant  \color{ekinorange} $x_0$}
        \State {Input parameters \color{ekinorange} $p_1,\ldots,p_k$}
        \State $ x_{1}  \gets  \Phi_1(; p_1,\ldots,p_k; x_0)$
        \State $ x_{2}  \gets  \Phi_2(x_1; p_1,\ldots,p_k; x_0)$
        \State \hspace{.3in} $\vdots$\hspace{.3in} $\vdots$
        \State $ x_{N}  \gets  \Phi_N(x_1,\ldots,x_{N-1}; p_1,\ldots,p_k; x_0)$

        \State {\color{ekinorange} Output $x_{N}$}  \end{algorithmic}
    \end{algorithm}
  \end{center}

  Algorithm \ref{gen} is an explicit computation.
  The {\it function} $\phi_i$ computes the
  value of the {\it variable} $x_i$.
  The notation
$\frac{d\Phi_i}{dx_j}$ or
$\frac{d\Phi_i}{dp_j}$
  gives the partial derivatives of one step
  of the algorithm.  By contrast, the
  notation $\frac{dx_i}{dp_j}$
  gives the partial derivatives across multiple
  steps of the algorithm.
  Algorithm \ref{gen}  is the general case of Algorithm~\ref{gen1}, the $\frac{d\Phi_i}{dx_j}$ and
$\frac{d\Phi_i}{dp_j}$ are general cases of what is seen
  in Algorithm~\ref{gen2}, and the $\frac{dx_i}{dp_j}$ generalize
  what is seen in Algorithm~\ref{gen3}.

  We note that the adjoint method literature tends to consider a yet more general implicit approach. Placing Section 3 of  \cite{johnson2007notes} in an explicit setting, we define a function $f$ such that $f(x, p)=0$.
  To this end, let
  \begin{equation}
    \label{explicit}
    f(x,p)= x - \Phi(x, p) := \begin{pmatrix} x_1 \\  x_2 \\ \vdots \\ x_N  \end{pmatrix}-
    \begin{pmatrix}
      \Phi_1(; p; x_0)    \\
      \Phi_2(x_1; p; x_0) \\
      \vdots              \\
      \Phi_N(x_1,\ldots,x_{N-1}; p; x_0)
    \end{pmatrix}
    .
  \end{equation}

  Clearly, given $p$, the computed  $x=(x_1,\ldots,x_N)$ from Algorithm \ref{gen} is a solution  to $f(x,p)=0$.
  Our goal is to reproduce \eqref{likeadjoint2}, which is the derivative of $x_N$ w.r.t. to the parameter $p$. 

  Let us first consider the derivation for $x_p$, which is the derivative of $x$, implicitly defined by $f(x, p) =  0$, w.r.t. to $p$. To connect the viewpoints a table of notation for various Jacobians is helpful:

\vspace{.2in}

\begin{center}
  \begin{tabular}{llll}
    \toprule
    Adjoint Method & Nabla Notation & Matrix                  & Size         \\
    \midrule
    $f_x$          & $\nabla_x f $  & $  I- \tilde{L} $       & $N\times N$  \\[.1in]
    $f_p$          & $\nabla_p f$   & $-M^T $                 & $N \times k$ \\[.1in]
    $x_p$          & $\nabla_p x$   & $(I-\tilde{L})^{-1}M^T$ & $N \times k$ \\
    \bottomrule
  \end{tabular}
\end{center}

\vspace{.1in}
The  matrices  themselves are explicitly:

\begin{displaymath}
  \tilde{L} = \left[\dfrac{\partial \Phi_i}{\partial x_j}\right]_{i,j}, \quad i > j,\ j = 1, \ldots, N-1,
\end{displaymath}
and

\[
  M^T = \left[\dfrac{\partial \Phi_i}{\partial p_j}\right]_{i,j}, \quad \nabla_p x =  \left[\dfrac{\partial x_i}{\partial p_j}\right]_{i,j}, \quad i \in 1,\ldots,N, j \in 1,\ldots,k.
\]

The matrix $\tilde{L}$ that contains the partials $\partial \Phi_j / \partial x_j$ is strictly lower triangular exactly because
Algorithm \ref{gen} is an explicit computation, whereas an implicit function would generally have a dense Jacobian.
Since $f(x,p) = x - \Phi(x,p)$, the Jacobian $\nabla_x f= I - \tilde{L}$.
Differentiating $0=f(x,p)$ with respect to $p$ we get $0 = f_x x_p + f_p$ or
$x_p = - f_x^{-1}f_p$ which is $(I-\tilde{L})^{-1}M^T$ in matrix notation explaining the bottom row of the above table.

If $g(x)$ is any scalar function of $x$, then the key adjoint equation is
\[
  \nabla_p g = g_xx_p = -g_x f_x^{-1} f_p := -\lambda^T f_p,
\]
where $\lambda$ satisfies the so-called adjoint equation $f_x^T \lambda =g_x^T$. Since $g_x$ is an $1$ by $k$ vector, by computing the adjoint $\lambda$ first, we reduce the computation of a matrix-matrix multiplication and a matrix-vector multiplication to two matrix-vector multiplications.

If we take $g(x)=x_N$ then $g_x = [0,\ldots,0,1]$.
The gradient is then
$$ \nabla_p g(x) =  [0,\ldots,0,1]  (I-\tilde{L})^{-1} M^T,$$
achieving our goal of reproducing \eqref{likeadjoint}. 

So much is happening here that it is worth repeating with other notation.
We can use the Jacobian of $f$ with respect to $x$ and $p$ to  differentiate \eqref{explicit}:
\[
  0=
  \begin{pmatrix}
    dx_1 \\ dx_2 \\  \vdots  \\ dx_N
  \end{pmatrix}
  -
  \begin{pmatrix}
    0                                     & 0      & \ldots                                     & 0      \\
    \dfrac{\partial \Phi_2}{\partial x_1} & 0      & \ldots                                     & 0      \\
    \vdots                                & \ddots & \vdots                                     & \vdots \\
    \dfrac{\partial \Phi_N}{\partial x_1} & \ldots & \dfrac{\partial \Phi_N}{\partial x_{N-1 }} & 0
  \end{pmatrix}
  \begin{pmatrix}
    dx_1 \\ dx_2 \\  \vdots \\ dx_N
  \end{pmatrix}
  -
  \begin{pmatrix}
    \dfrac{\partial \Phi_1}{\partial p_1} & \ldots & \dfrac{\partial \Phi_1}{\partial p_k } \\
    \vdots                                & \vdots & \vdots                                 \\
    \dfrac{\partial \Phi_N}{\partial p_1} & \ldots & \dfrac{\partial \Phi_N}{\partial p_k }
  \end{pmatrix}
  \begin{pmatrix}
    dp_1 \\ dp_2 \\  \vdots \\ dp_k
  \end{pmatrix},
\]
which can be solved to obtain
\[
  \begin{pmatrix}
    dx_1 \\ dx_2 \\  \vdots \\ dx_N
  \end{pmatrix}
  =
  \left(
  I -
  \begin{pmatrix}
    0                                     & 0      & \ldots                                     & 0      \\
    \dfrac{\partial \Phi_2}{\partial x_1} & 0      & \ldots                                     & 0      \\
    \vdots                                & \ddots & \vdots                                     & \vdots \\
    \dfrac{\partial \Phi_N}{\partial x_1} & \ldots & \dfrac{\partial \Phi_N}{\partial x_{N-1 }} & 0
  \end{pmatrix}
  \right)^{-1}
  \begin{pmatrix}
    \dfrac{\partial \Phi_1}{\partial p_1} & \ldots & \dfrac{\partial \Phi_1}{\partial p_k } \\
    \vdots                                & \vdots & \vdots                                 \\
    \dfrac{\partial \Phi_N}{\partial p_1} & \ldots & \dfrac{\partial \Phi_N}{\partial p_k }
  \end{pmatrix}
  \begin{pmatrix}
    dp_1 \\ dp_2 \\  \vdots \\ dp_k
  \end{pmatrix} .
\]

Some readers unfamiliar with the notation of differentials might prefer what amounts to a notational change, but avoids
the notation of differentials:
\begin{displaymath}
  \begin{pmatrix}
    \dfrac{\partial x_1}{\partial p_1} & \ldots & \dfrac{\partial x_1}{\partial p_k} \\
    \vdots                             & \vdots & \vdots                             \\
    \dfrac{\partial x_N}{\partial p_1} & \ldots & \dfrac{\partial x_N}{\partial p_k} \\
  \end{pmatrix}
  =
  \left(
  I -
  \begin{pmatrix}
    0                                      & 0      & \ldots                                    & 0      \\
    \dfrac{\partial \Phi_2}{\partial x_1} & 0      & \ldots                                    & 0      \\
    \vdots                                 & \ddots & \vdots                                    & \vdots \\
    \dfrac{\partial \Phi_N}{\partial x_1} & \ldots & \dfrac{\partial \Phi_N}{\partial x_{N-1}} & 0
  \end{pmatrix}
  \right)^{-1}
  \begin{pmatrix}
    \dfrac{\partial \Phi_1}{\partial p_1} & \ldots & \dfrac{\partial \Phi_1}{\partial p_k } \\
    \vdots                                 & \vdots & \vdots                                 \\
    \dfrac{\partial \Phi_N}{\partial p_1} & \ldots & \dfrac{\partial \Phi_N}{\partial p_k }
  \end{pmatrix} .
\end{displaymath}


\section{Julia, the power of language}

\subsection{The challenge}
\label{challenge}

This section provides a complete realization of the challenge described in the preface (Section \ref{preface}).
The question we asked is whether we
could  bring to life the linear algebra mathematics
expressed in
$$\nabla J = M^T((I-L)^T \backslash g)$$
by typing
the command
\begin{center}
  \includegraphics[width=1.5in]{figures/goal.png} 
\end{center}
and  computing the
backpropagated gradient of a matrix neural network almost by magic?

We remark that it is common to see code  in papers.
Code can serve the purpose of specifying details, facilitating
reproducibility, and verifiability.  Code can also allow users
to adapt methods to their own situations.
In addition to all of the above, we have a further purpose.
We believe the code example we provide shows
the power, elegance, and utility of the Julia programming language
in ways that may be difficult or impossible to imagine in
other languages. 

At the risk of showing the end of the code before the start,
63  lines of setup culminate in exactly what we wanted:
code which looks just like the math of matrices with operators
that correctly calculates the gradient fulfilling our
title goal of backpropagating through back substitution with a backslash:

\begin{center}
  \includegraphics[width=\textwidth]{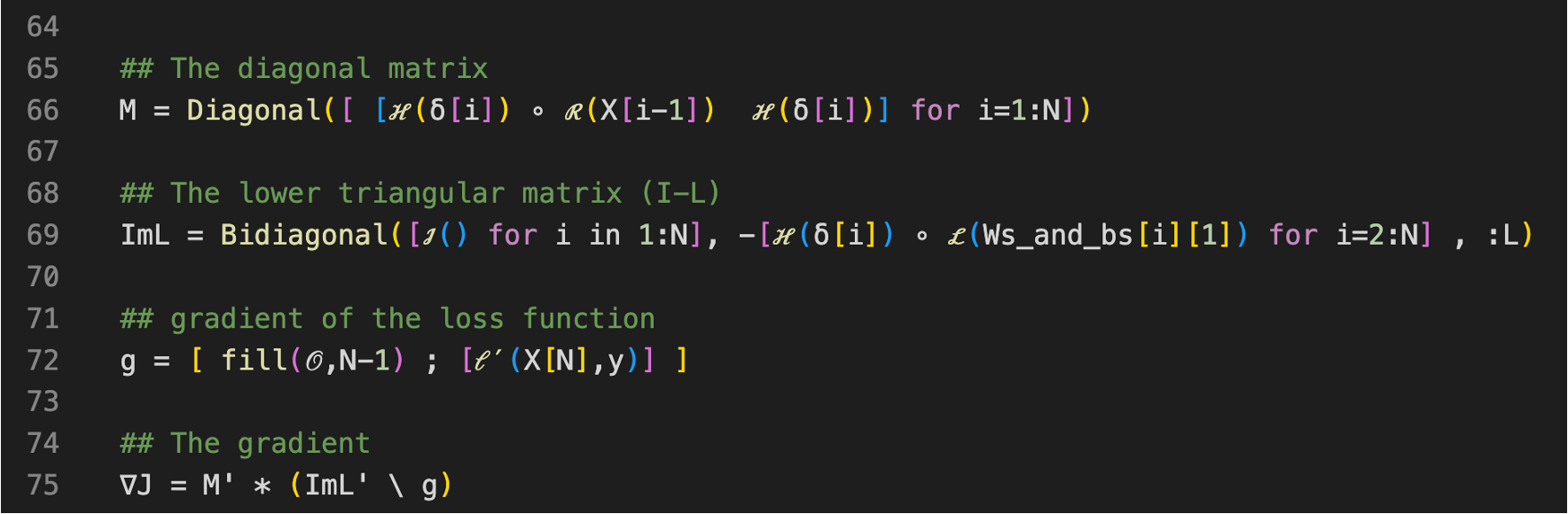}
\end{center}

The first 28 lines elegantly set up
the mathematics very much like a mathematician defining
operators and algebraic axioms:

\begin{center}
  \includegraphics[width=0.92\textwidth]{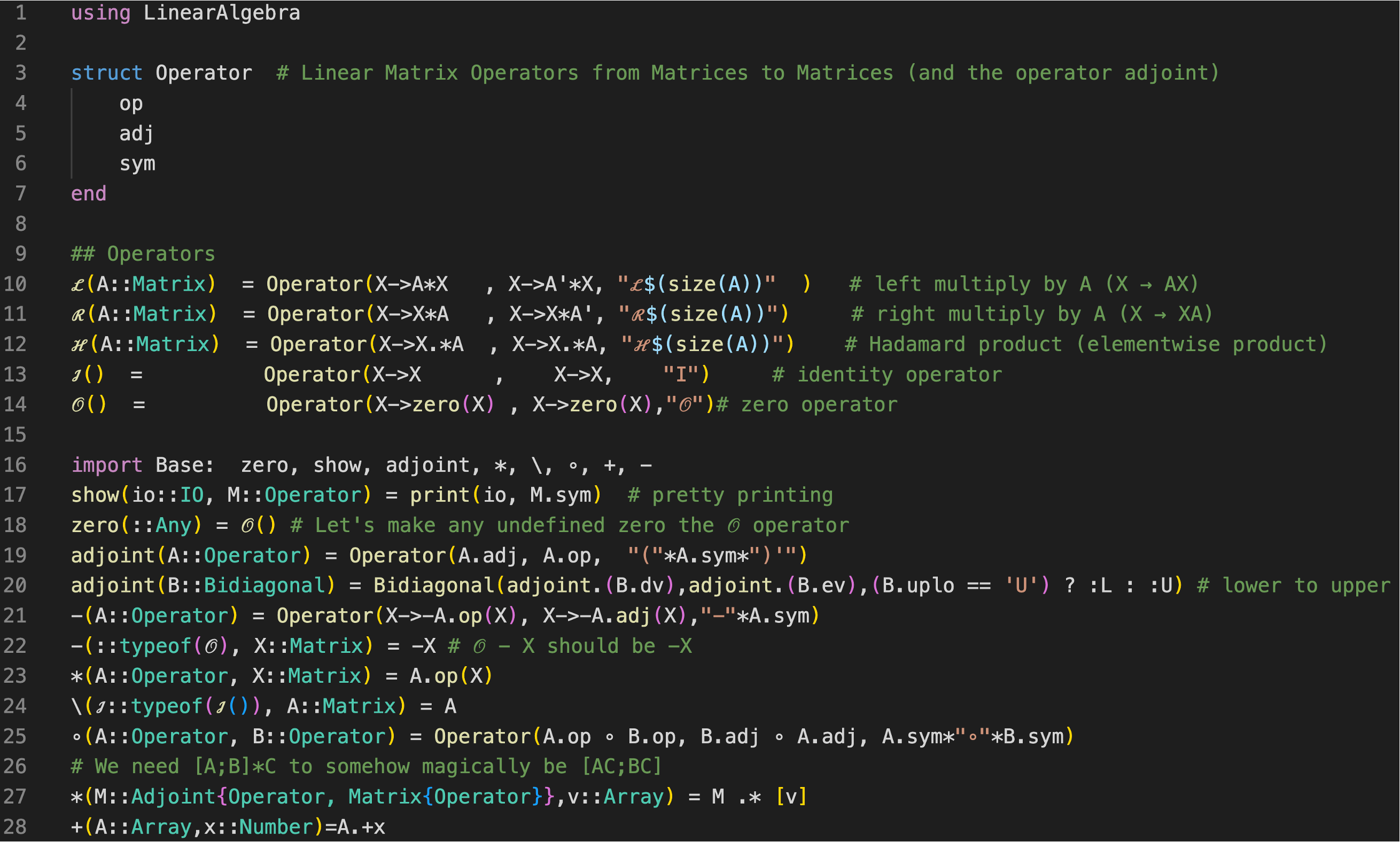}
\end{center}

Lines 10-14  above define matrix operators and their adjoints.
Lines 16-28 define various math operations, such
as the negative operator on line 21, or the composition of operators
on line 25.
\begin{center}
\includegraphics[width=.75\textwidth]{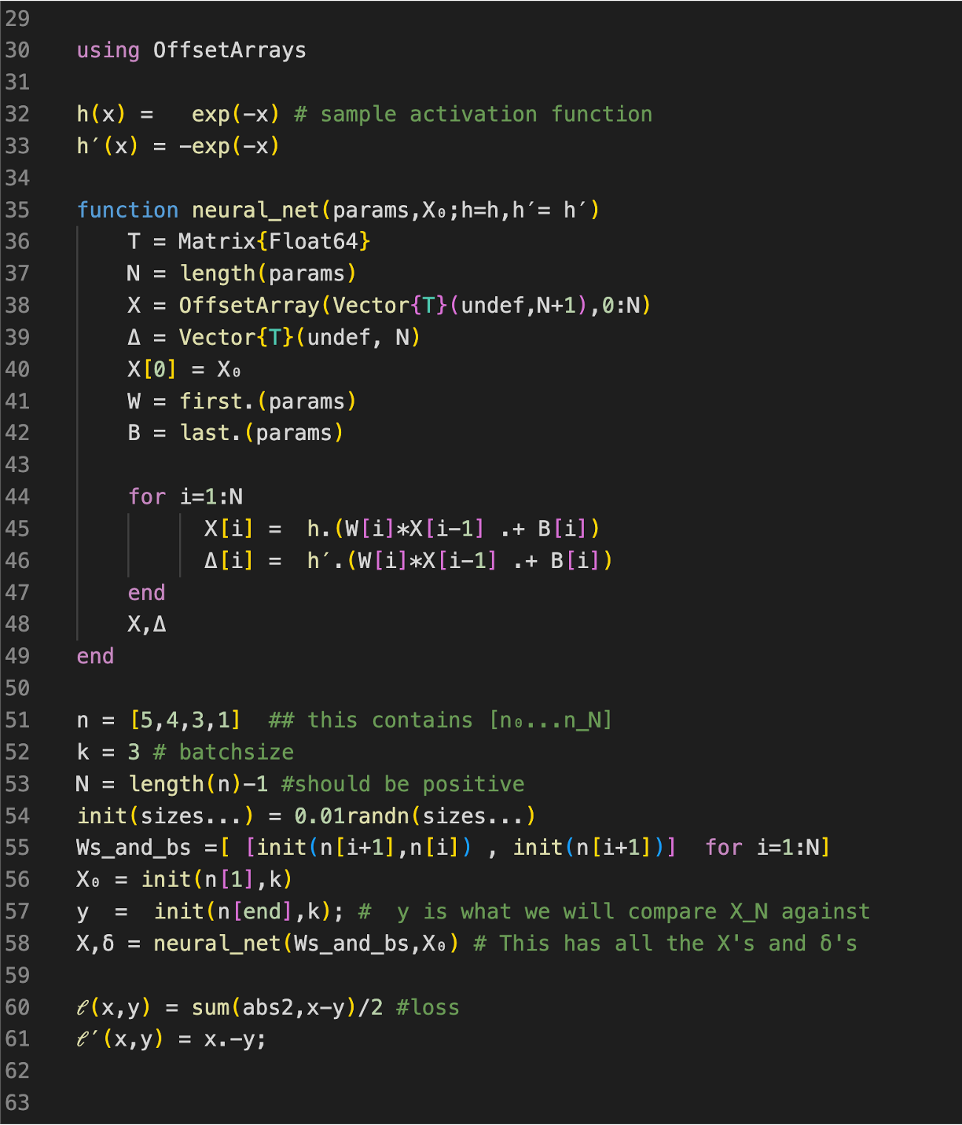}
\end{center}
For completeness  we list   lines 29 through 63 which
constitute the setup of
a basic forward pass through a matrix neural  net.  We remark  that lines 30 and  38 allow
an index origin of 0. The readers are encouraged to try the code at \small{\url{https://github.com/alanedelman/BackBackBack}}.

\subsection{Modern Computer Science meets Linear Algebra}

The venerable position of numerical linear algebra libraries
can not be undersold.  Years of rigorous mathematical and
algorithmic research have culminated in the modern LAPACK library
\cite{anderson1999lapack} which represents
a crowning achievement of value to a huge number
of users who call LAPACK perhaps from, for example, Julia, NumPy, or MATLAB.
In most cases the users are unaware of the scientific bedrock of which they
are beneficiaries.

Continuing this grand tradition, we wish to highlight some
of the computer science innovations that allow for the code
in Section \ref{challenge} to look so deceptively simple.\\

\noindent\emph{Generic Programming or how can the backslash just work? }
We invite the reader to consider how the innocent backslash on line 75
of the code in Section \ref{challenge} could possibly perform a backpropogation of derivative.
We believe this would be impossible in, for example, NumPy or MATLAB
as these packages currently exist.
From a computer science point of view, Julia's multiple dispatch
mechanism and generic programming features allow the
generic backslash to work with matrices and vectors whose elements
are operators and compositions of operators.  We remind the reader
that the operators are not language constructs, but are created
in software on the first 28 lines of code.  The backslash,
however, is not LAPACK's backslash, as the LAPACK library
is constrained to floating point real and complex numbers.
Julia's backslash currently runs LAPACK when dispatched by
matrices of floats, but, as is the case here, the generic algorithm
is called.
We are fascinated by the fact that the author of the generic algorithm
would not have imagined how it might be used.  We are aware of
backslash being run on quaternion matrices, block matrices, matrices over finite fields,  and now matrices with operators.  Such is the mathematical power
of abstraction and what becomes possible if software is allowed to be generic.
In the context of backpropagation, replacing the ``for loops'' with
the backslash helps us see backpropogation from a higher viewpoint.\\

\noindent\emph{The significance of transpose all the way down.}
Not without controversy, Julia implements transpose recursively.  We believe this
is the preferred behavior.  This means a block matrix of block matrices of matrices (etc.) will transpose in the expected manner.  Similarly matrices
of complex number or quaternions will perform conjugate transposes as expected.
In this work the $M$ as seen in Line 66 of the code in Section
\ref{challenge} is diagonal, but is not symmetric.
In line 75 we are transposing a diagonal matrix of $1\times 2$ matrices
of composed operators \verb+M'+  while
in that same line we are also transposing a bidiagonal matrix  of  operators.
Because the operator adjoint is defined on lines 10-14 of the code
and the adjoint for a composed operator is defined on line 25,
Julia's generic implementation, again, just works.  We are not
aware of any other linear algebra system whereby the transpose
would just work this readily.
The page {{\small{\noindent
      \verb+https://discourse.julialang.org/t/why-is-transpose-recursive/2550+}}} documents some of the controversy. We are extremely grateful that the recursive
definition won the day.\\

\noindent\emph{A quick word about performance.}
There is nothing in the backslash formulation that would impede performance.\\

\noindent\emph{Possible extensions to  the example code  in Section \ref{challenge}.}
We deliberately only used as an example the matrix neural network.
We also have implemented a fully connected neural network where the
matrix $I-L$ is  a Julia triangular type, whereas the reference example
was bidiagonal.   We also implemented a square case where the $W$ parameter
was constant from one iteration to the next.
We also conceived of the case of being restricted to a manifold.
We thus stress that
we did not build a fully functional package at this time, and
thus emphasize that this could be future research, but we have not
yet seen any roadblock to this methodology.\\

\noindent\emph{Concluding Moral.}
Exciting new innovations in numerical algorithms are emerging from software developments. Critical
elements 
for creativity 
include: generic programming (generic operators),
abstract representations, 
fast performance without waste, 
multiple dispatch, and an aggressive type system. 

Abstraction matters. Software matters. Language matters.

\section{Acknowledgments} 
\small{We would like to thank the anonymous reviewers for their valuable comments. We wish to thank David Sanders and Jeremy Kepner for  helpful conversations. This material is based upon work supported by the National Science Foundation under grant no. OAC-1835443, grant no. SII-2029670, grant no. ECCS-2029670, grant no. OAC-2103804, and grant no. PHY-2021825. We also gratefully acknowledge the U.S. Agency for International Development through Penn State for grant no. S002283-USAID. The information, data, or work presented herein was funded in part by the Advanced Research Projects Agency-Energy (ARPA-E), U.S. Department of Energy, under Award Number DE-AR0001211 and DE-AR0001222. This material is based upon work supported by the Defense Advanced Research Projects Agency (DARPA) under Agreement No HR00112290091. We also gratefully acknowledge the U.S. Agency for International Development through Penn State for grant no. S002283-USAID. The views and opinions of authors expressed herein do not necessarily state or reflect those of the United States Government or any agency thereof. This material was supported by The Research Council of Norway and Equinor ASA through Research Council project "308817 - Digital wells for optimal production and drainage". Research was sponsored by the United States Air Force Research Laboratory and the United States Air Force Artificial Intelligence Accelerator and was accomplished under Cooperative Agreement Number FA8750-19-2-1000. The views and conclusions contained in this document are those of the authors and should not be interpreted as representing the official policies, either expressed or implied, of the United States Air Force or the U.S. Government. The U.S. Government is authorized to reproduce and distribute reprints for Government purposes notwithstanding any copyright notation herein.}

\bibliographystyle{siamplain}
\bibliography{references}

\begin{thebibliography}{10}

\bibitem{anderson1999lapack}
{\sc E.~Anderson, Z.~Bai, C.~Bischof, L.~S. Blackford, J.~Demmel, J.~Dongarra, J.~Du~Croz, A.~Greenbaum, S.~Hammarling, A.~McKenney, et~al.}, {\em LAPACK users' guide}, SIAM, 1999.

\bibitem{JMLR:v18:17-468}
{\sc A.~G. Baydin, B.~A. Pearlmutter, A.~A. Radul, and J.~M. Siskind}, {\em Automatic differentiation in machine learning: a survey}, Journal of Machine Learning Research, 18 (2018), pp.~1--43, \url{http://jmlr.org/papers/v18/17-468.html}.

\bibitem{bezanson2017julia}
{\sc J.~Bezanson, A.~Edelman, S.~Karpinski, and V.~B. Shah}, {\em Julia: A fresh approach to numerical computing}, SIAM review, 59 (2017), pp.~65--98.

\bibitem{bradley2010pde}
{\sc A.~M. Bradley}, {\em {PDE}-constrained optimization and the adjoint method}, 2010, \url{https://cs.stanford.edu/~ambrad/adjoint_tutorial.pdf}.

\bibitem{Goodfellow-et-al-2016}
{\sc I.~Goodfellow, Y.~Bengio, and A.~Courville}, {\em Deep Learning}, MIT Press, 2016.
\newblock \url{http://www.deeplearningbook.org}.

\bibitem{griewank2003mathematical}
{\sc A.~Griewank}, {\em A mathematical view of automatic differentiation}, Acta Numerica, 12 (2003), pp.~321--398.

\bibitem{griewank2008evaluating}
{\sc A.~Griewank and A.~Walther}, {\em Evaluating derivatives: principles and techniques of algorithmic differentiation}, SIAM, 2008.

\bibitem{Wik:DifProg}
{\sc M.~Innes, A.~Edelman, K.~Fischer, C.~Rackauckus, E.~Saba, V.~Shah, and W.~Tebbutt}, {\em $\partial{P}$: A differentiable programming system to bridge machine learning and scientific computing}, 2019, \url{http://arxiv.org/abs/1907.07587}.

\bibitem{johnson2007notes}
{\sc S.~G. Johnson}, {\em Notes on adjoint methods for 18.335}, 2006, \url{https://math.mit.edu/~stevenj/18.336/adjoint.pdf}.

\bibitem{kepner2011graph}
{\sc J.~Kepner and J.~Gilbert}, {\em Graph algorithms in the language of linear algebra}, SIAM, 2011.

\bibitem{li2018differentiable}
{\sc T.-M. Li, M.~Gharbi, A.~Adams, F.~Durand, and J.~Ragan-Kelley}, {\em Differentiable programming for image processing and deep learning in halide}, ACM Transactions on Graphics (TOG), 37 (2018), pp.~1--13.

\bibitem{1810.08297}
{\sc J.~Revels, T.~Besard, V.~Churavy, B.~D. Sutter, and J.~P. Vielma}, {\em Dynamic automatic differentiation of {GPU} broadcast kernels}, 2018, \url{https://arxiv.org/abs/arXiv:1810.08297}.

\bibitem{revels2016forward}
{\sc J.~Revels, M.~Lubin, and T.~Papamarkou}, {\em Forward-mode automatic differentiation in {J}ulia}, arXiv preprint arXiv:1607.07892,  (2016).

\bibitem{stanley2015catalan}
{\sc R.~P. Stanley}, {\em Catalan numbers}, Cambridge University Press, 2015.

\end{thebibliography}
\end{document}